\documentclass[10pt]{article}
\usepackage[a4paper,margin=29mm]{geometry}
\usepackage{amsmath,amssymb,amsthm,mathtools}
\usepackage{graphicx,booktabs,array}
\usepackage{tikz}
\usepackage{microtype}
\usepackage[shortlabels]{enumitem}
\usepackage[colorlinks=true,allcolors=black]{hyperref}
\usepackage[nameinlink,noabbrev]{cleveref}

\makeatletter
\renewcommand{\@maketitle}{%
  \null\vskip 1em
  \begin{center}
    {\normalfont\normalsize\bfseries\MakeUppercase{\@title}\par}
    \vskip 1.2em
    {\normalfont\footnotesize\MakeUppercase{\@author}\par}
  \end{center}
  \vskip .6em
}
\renewcommand{\@seccntformat}[1]{\textup{\csname the#1\endcsname.}\quad}
\renewcommand{\section}{\@startsection{section}{1}{\z@}%
  {-3.5ex \@plus -1ex \@minus -.2ex}{1.5ex \@plus .2ex}%
  {\normalfont\normalsize\scshape\centering}}
\renewcommand{\subsection}{\@startsection{subsection}{2}{\z@}%
  {-3.25ex \@plus -1ex \@minus -.2ex}{1ex \@plus .2ex}%
  {\normalfont\normalsize\bfseries}}
\renewcommand{\subsubsection}{\@startsection{subsubsection}{3}{\z@}%
  {-3.25ex \@plus -1ex \@minus -.2ex}{1ex \@plus .2ex}%
  {\normalfont\normalsize\itshape}}
\makeatother
\renewenvironment{abstract}
  {\footnotesize
   \list{}{\leftmargin=2.5em\rightmargin=2.5em
     \topsep=0pt\partopsep=0pt\parsep=0pt\itemsep=0pt}
   \item\relax\noindent\textsc{\abstractname.}\enspace\ignorespaces}
  {\endlist\vspace{1em}}

\allowdisplaybreaks[1]
\newtheorem{maintheorem}{Theorem}

\crefname{maintheorem}{Theorem}{Theorems}
\Crefname{maintheorem}{Theorem}{Theorems}
\newtheorem*{introtheoremA}{Theorem~\ref{thm:masked-maclaurin}}
\newtheorem*{introtheoremB}{Theorem~\ref{sg:theorem}}
\newtheorem{theorem}{Theorem}[section]
\newtheorem{proposition}[theorem]{Proposition}
\newtheorem{lemma}[theorem]{Lemma}
\newtheoremstyle{starredplain}{\topsep}{\topsep}
  {\itshape}{}{\bfseries}{.}{5pt plus 1pt minus 1pt}
  {\thmname{#1}\thmnumber{ *#2}\thmnote{ (#3)}}
\theoremstyle{starredplain}
\newtheorem{starredlemma}[theorem]{Lemma}
\crefname{starredlemma}{Lemma}{Lemmas}
\Crefname{starredlemma}{Lemma}{Lemmas}
\theoremstyle{plain}
\newtheorem{corollary}[theorem]{Corollary}
\theoremstyle{definition}
\newtheorem{definition}[theorem]{Definition}
\theoremstyle{remark}
\newtheorem{remark}[theorem]{Remark}
\newcommand{\R}{\mathbb R}

\newcommand{\Cl}{\operatorname{Cl}}
\newcommand{\Pf}{\operatorname{Pf}}
\newcommand{\rank}{\operatorname{rank}}

\newcommand{\Kcal}{\mathcal K}

\title{Discrete Helicity on Graphs: \\ Pfaffian Energies and Clique Bounds}
\author{Jingqi Li}
\date{}
\begin{document}
\maketitle
\begin{abstract}
We introduce a helicity functional of each order on graphs for real-valued edge forms, and give it a unified local Pfaffian realization. Motivated by the continuous energy inequality, we further construct a discrete version of helicity energy, and prove sharp Maclaurin inequalities for these energies. In particular, a special case of these Maclaurin inequalities yields a lower bound on the clique number, depending on the choice of edge form. As an application, we use this certificate to prove that every nonempty $d$-regular graph with smallest adjacency eigenvalue at least $-3$ satisfies $\omega(G)\ge d\,g(\eta)-C$, where $C$ is an absolute constant and $g(\eta)$ is a rank-dependent coefficient taking values between $1/4$ and $1/3$. Building on the structural results of Koolen, Yang and Yang~\cite{KYY}, this improves the known clique-number bound to the optimal leading constant $1/4$, whose asymptotic sharpness follows from Yang and Koolen's construction~\cite{YangKoolen2025}.

\end{abstract}

\setcounter{tocdepth}{2}
\begingroup
\hypersetup{linkcolor=blue!60!black,linktoc=all}
\tableofcontents
\endgroup
\clearpage

\section{Introduction}
\label{sec:introduction}
Helicity originated in the study of the linkage of vortex lines and
magnetic field lines, and has become a fundamental invariant in
topological hydrodynamics. On a closed oriented three-manifold
equipped with a volume form, helicity is defined for exact
divergence-free vector fields and is invariant under
volume-preserving diffeomorphisms
\cite{Moffatt1969,Arnold1986,ArnoldKhesin2021}.
It therefore provides a way to extract global topological
information from a field whose behaviour is described locally
by differential equations.

The significance of helicity is reflected in its connections with
energy, dynamics and geometry. Once a Riemannian metric is fixed,
Arnold's energy inequality bounds the energy of a field from below
in terms of its helicity \cite{Arnold1986}. Nonzero helicity consequently obstructs
relaxation to a state of arbitrarily small energy. On a homology
three-sphere, Arnold's helicity theorem identifies helicity with
the averaged asymptotic linking of trajectories, relating an
integral invariant to the long-time behaviour of the flow
\cite{Arnold1986,ArnoldKhesin2021}.
Its differential-form expression also coincides, up to
normalization, with the Abelian Chern--Simons functional
\cite{ChernSimons1974}.
Related helicity-type integrals in higher odd dimensions involve
densities of the form $\alpha\wedge(d\alpha)^{s-1}$ and admit
geometric interpretations through the linking of foliations
\cite{Khesin1992}.
These connections suggest that the algebra underlying helicity
may remain useful beyond its original hydrodynamical setting.

The discretization of helicity has been studied extensively in
computational electromagnetics and numerical partial differential
equations. Early work of Kotiuga investigated helicity functionals,
their metric invariance and the finite-element matrices arising
from their discretization \cite{Kotiuga1989,Kotiuga1990}.
Discrete exterior calculus (DEC), developed by Hirani and others
\cite{Hirani2003,DesbrunHiraniLeokMarsden2005}, and finite element
exterior calculus (FEEC) \cite{ArnoldFalkWinther2006} subsequently
provided systematic frameworks for retaining differential-form
and cohomological structures under discretization.
Within finite element de Rham complexes, Hu, Lee and Xu constructed
schemes for incompressible magnetohydrodynamics that preserve
magnetic and cross helicity in the ideal, unforced setting under
appropriate boundary conditions \cite{HuLeeXu2021}.
More recently, He, Farrell, Hu and Andrews developed a
discretization for magnetic relaxation that preserves helicity
while reproducing energy dissipation, and established a discrete
Arnold inequality \cite{HeFarrellHuAndrews2026}.
Such results demonstrate the practical importance of retaining
the relationship between helicity and energy at the discrete level.

The numerical constructions just described approximate fields on
an underlying geometric domain. The present work arose from
asking how much of the helicity construction survives when that
geometric setting is absent. Our starting observation is local:
on an ordered tetrahedron, the integrated Whitney helicity density
is a fixed multiple of the Pfaffian of the skew-symmetric matrix
of its six edge values. This expression depends only on the edge
data and the vertex ordering. The observation led us to ask whether
helicity and its higher-dimensional analogues could be formulated
directly on an abstract graph, and whether their local algebra
could reveal properties of the graph itself.

A finite graph already carries the requisite simplicial structure
through its clique complex. Cochains, coboundaries and their
products can therefore be considered without assuming that the
complex triangulates a manifold. We use this structure to define
helicity functionals by evaluation on specified simplicial cycles,
and study their local expressions through Whitney forms.
Our first result is a unified Pfaffian formula for these local
helicities in every order. For every finite simple graph $G$,
integer $s\ge1$, $2s$-vertex clique $S$ with its reference
orientation, and real edge form $a\in C^1(G)$,
\[
\mathcal H_{s,S}(a)
=\frac{2^{s-1}(s-1)!}{(2s-1)!}\Pf\bigl(a[S]\bigr).
\]
Here $\mathcal H_{s,S}(a)$ is the local helicity functional,
$a$ assigns real values $a_{ij}=-a_{ji}$ to oriented edges, and
$a[S]$ is their skew-symmetric matrix in the same vertex order.
Let $\Kcal_{2s}(G)$ denote the set of $2s$-vertex cliques of $G$.
Summing the squares of the local helicities over these cliques
yields a nonnegative discrete helicity energy. This construction
is motivated by the idea of approximating continuous energy
through local data on fine triangulations of manifolds.
Dividing this energy by the square of the displayed coefficient
gives the normalized helicity energy $T_{s,G}(a)$ in
\eqref{eq:energy-definition}; equivalently, one sums the squares
of the normalized local helicities, which are exactly Pfaffians.
Our principal interest is in the relationship between these local
helicity densities, the energy of an edge form and the combinatorial
structure of the graph.

The two main results in our paper develop this relationship quantitatively.

\begin{introtheoremA}
Let $G$ be a finite simple graph and put $p=\lfloor\omega(G)/2\rfloor$. For every edge form $a\in C^1(G)$ and all integers $1\le j<k\le p$,
\[
\left(\frac{T_{k,G}(a)}{\binom{p}{k}}\right)^{1/k}
\le
\left(\frac{T_{j,G}(a)}{\binom{p}{j}}\right)^{1/j}.
\]
The inequality is sharp.
\end{introtheoremA}

Its specialization to $j=1$ and $k=2$ yields an edge-form certificate for the existence of a large clique; see \eqref{eq:certificate-introduction}. Combining this certificate with the structural theory of Koolen, Yang and Yang \cite{KYY}, we obtain the following rank-dependent bound.

\begin{samepage}
\begin{introtheoremB}
There is an absolute constant $C_*$ such that every nonempty finite simple $d$-regular graph $G$ on $n$ vertices, with adjacency matrix $M_G$ and $\lambda_{\min}(M_G)\ge-3$, satisfies
\[
\omega(G)\ge d\,g(\eta)-C_*,\qquad
\eta=\frac{\rank(M_G+3I)}{n},
\]
where
\[
g(\eta)=
\begin{cases}
\dfrac{1-2\eta}{3-4\eta},&0\le\eta\le\tfrac14,\\[4pt]
\dfrac{2\eta+1}{6},&\tfrac14\le\eta\le\tfrac12,\\[4pt]
\dfrac13,&\tfrac12\le\eta\le1.
\end{cases}
\]
\end{introtheoremB}
\end{samepage}

The coefficient $g(\eta)$ ranges from $1/4$ to $1/3$, so in particular $\omega(G)\ge d/4-C_*$. The construction of Yang and Koolen \cite{YangKoolen2025} shows that the universal leading constant $1/4$ is optimal. The full statement in \cref{sg:section} also supplies an edge form certifying the bound for all sufficiently large $d$.

The paper is organized as follows.
\Cref{sec:preliminaries} introduces the basic definitions and recalls
the classical results used throughout the paper.
\Cref{sec:helicity} develops helicity on clique complexes and
establishes its local Pfaffian interpretation, including the
higher-dimensional formula.
\Cref{sec:maclaurin} proves the sharp helicity energy Maclaurin inequality.
\Cref{sg:section} applies the case $j=1$, $k=2$ to obtain the
rank-dependent clique bound.
Finally, \cref{sec:questions} discusses further questions concerning
helicity, energy and graph structure.

\begin{remark}
A broader motivation is the possibility that helicity may supply
useful tools for combinatorial problems beyond those considered
here. The entropy method offers an instructive precedent:
an invariant arising in another subject can become an effective
instrument for proving combinatorial inequalities, as illustrated
by Radhakrishnan's entropy proof of Bregman's theorem
\cite{Radhakrishnan1997}.
We hope that the local algebra of helicity can play a similar
role. The passage from a tetrahedral helicity identity to a
spectral clique bound provides one concrete instance of this
possibility.

Although helicity originates in physics and topological hydrodynamics, this paper is essentially combinatorial and requires no topological machinery beyond a limited use of Whitney interpolation on abstract simplicial complexes. However, in order to accommodate readers with different backgrounds, we nevertheless provide a detailed introduction in \cref{sec:preliminaries} to both the graph-theoretic and the differential-form preliminaries.
\end{remark}

\begin{remark}
Graphs $G$ with adjacency matrix $M_G$ and $\lambda_{\min}(M_G)\ge-3$
have been studied through structural descriptions, integral representations
and clique bounds. Koolen, Yang and Yang~\cite{KYY} developed a
Hoffman-graph description and proved $2$-integrability, an integral
Gram representability property after scaling, for connected graphs
in this class with sufficiently large minimum degree. Under the
additional assumption that every pair of vertices at distance two
has the same number $c\ge9$ of common neighbors, Yang et
al.~\cite{YangEtAl2026} proved $1$-integrability and classified the
connected regular graphs with $\lambda_{\min}(M_G)\in[-3,-2)$
and sufficiently large degree. On the extremal side, Yang and
Koolen~\cite[Corollary~1.5]{YangKoolen2025} showed that, for sufficiently
large $d$, every vertex of a connected $d$-regular graph with
$\lambda_{\min}(M_G)\ge-3$ lies in a clique of order at least
$(d-\varepsilon)/15$, where $\varepsilon>0$ is an absolute constant.

This spectral condition also connects the graph class with lattice
theory and discrete geometry. Indeed, $M_G+3I$ is a positive
semidefinite integral Gram matrix whose vectors generate an integral
lattice~\cite{KYY}. Dividing these vectors by $\sqrt{3}$ gives unit
vectors with pairwise inner products in $\{0,1/3\}$, hence a spherical
set with at most two distances. Their span has dimension
$\rank(M_G+3I)$, giving a geometric meaning to the rank parameter
in our bound. Combinatorial design theory supplies regular examples
in this class, including block graphs of Steiner triple systems
and Latin square graphs~\cite{YangEtAl2026}. Steiner systems
$S(2,4,v)$ also underlie the extremal construction of Yang and
Koolen~\cite{YangKoolen2025}.
\end{remark}

\section{Preliminaries}
\label{sec:preliminaries}

This section collects the definitions, conventions and classical
results used in the paper. We begin with graphs, clique complexes
and spectral tools, then introduce simplicial chains and cochains.
We next review exterior algebra, Whitney interpolation and cochain
products, together with the classical local helicity formulas that
motivate our construction. Finally, we recall Pfaffians and the
Maclaurin inequalities used to compare the graph energies.

Readers familiar with graph theory and graph spectra may skip
\cref{subsec:graphs-clique-complexes}, while those familiar with
simplicial chains and cochains may skip
\cref{subsec:chains-cochains} after noting our orientation and
coboundary conventions. Readers comfortable with exterior algebra
and differential forms may begin their reading of
\cref{subsec:products-whitney} with Whitney interpolation.
 References for
further reading are provided alongside the relevant material in
this section and in the main text.

\subsection{Graphs, clique complexes, and spectra}
\label{subsec:graphs-clique-complexes}
\label{subsec:spectral-preliminaries}

The material in this subsection falls into three groups.
We first introduce the basic graph definitions and the clique
complex. We then recall the spectral tools needed later:
Rayleigh quotients, eigenvalue interlacing, Perron--Frobenius
theory, Smith's classification and the graph Laplacian.
Finally, we record the elementary consequence of Ramsey's
theorem used to find cliques inside vertex neighborhoods.

\subsubsection{Graphs and clique complexes}

Throughout this paper, $G=(V,E)$ is a finite simple undirected graph
with nonempty vertex set. We fix a total ordering of $V$ and identify
$V$ with $[n]=\{1,\ldots,n\}$ in this order. Each edge is an unordered
pair of distinct vertices; loops and multiple edges are excluded.
For $S\subseteq V$, the induced subgraph $G[S]$ has vertex set $S$
and contains precisely those edges of $G$ whose endpoints both lie
in $S$. The degree of a vertex is the number of its neighbours, and
$G$ is \emph{$d$-regular} if every vertex has degree $d$.
We denote the adjacency matrix of $G$ by $M_G$, with
$(M_G)_{ij}=1$ when $\{i,j\}\in E$ and $(M_G)_{ij}=0$ otherwise.
Thus $M_G$ is a real symmetric matrix with zero diagonal.

A \emph{clique} or \emph{complete subgraph} is a subset of vertices whose distinct elements are
pairwise adjacent. We write $K_r$ for the complete graph on $r$
vertices and $\mathcal K_r(G)$ for the set of all
$r$-cliques of $G$. In particular, $\mathcal K_1(G)$ consists of the
singleton vertex sets, and $\mathcal K_2(G)$ is the edge set.
The \emph{clique number} $\omega(G)$ is the largest size of a clique:
\[
  \omega(G)=\max\{r:\mathcal K_r(G)\neq\varnothing\}.
\]

An \emph{abstract simplicial complex} is a family of finite vertex
sets closed under taking subsets. Its elements are called
\emph{simplices}, and a nonempty simplex with $r$ vertices has
dimension $r-1$. Every graph determines such a complex, its
\emph{clique complex},
\[
  \operatorname{Cl}(G)
  =
  \bigl\{
    S\subseteq V:
    \{i,j\}\in E
    \text{ for all distinct } i,j\in S
  \bigr\}.
\]
The empty set is included in this definition. Closure under taking
subsets follows because every subset of a clique is again a clique.
Thus the $(r-1)$-dimensional simplices of $\operatorname{Cl}(G)$
correspond exactly to the members of $\mathcal K_r(G)$, and
\[
  \dim\operatorname{Cl}(G)=\omega(G)-1.
\]
The vertices and edges of this complex, collectively called its
\emph{$1$-skeleton}, recover the original graph $G$.

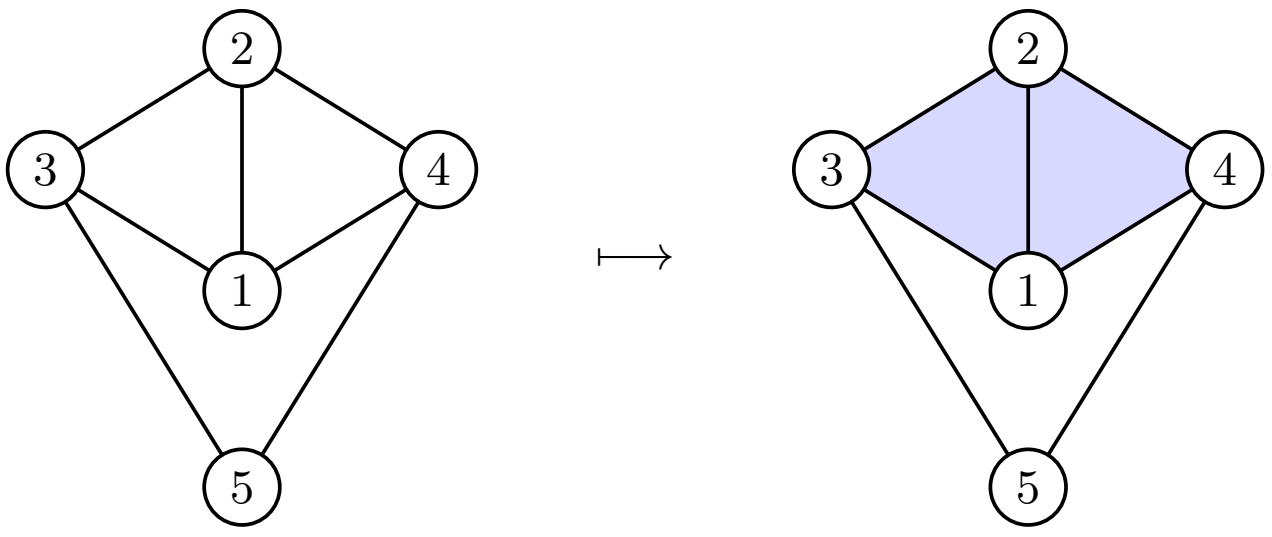
\begin{figure}[htbp]
  \centering
  \begin{tikzpicture}[
    line width=0.7pt,
    line cap=round,
    line join=round,
    font=\small,
    vertex/.style={
      circle,
      draw=black,
      fill=white,
      minimum size=5mm,
      inner sep=0pt
    }
  ]

    \begin{scope}
      \coordinate (g1) at (0,0);
      \coordinate (g2) at (0,1.6);
      \coordinate (g3) at (-1.3,0.8);
      \coordinate (g4) at (1.3,0.8);
      \coordinate (g5) at (0,-1.3);

      \foreach \i/\j in {1/2,1/3,2/3,1/4,2/4,3/5,4/5}
        \draw (g\i) -- (g\j);

      \foreach \i in {1,...,5}
        \node[vertex] at (g\i) {$\i$};

      \node at (0,-1.95) {(a) Graph $G$};
    \end{scope}

    \node at (2.6,0.2) {$\longmapsto$};

    \begin{scope}[xshift=5.2cm]
      \coordinate (k1) at (0,0);
      \coordinate (k2) at (0,1.6);
      \coordinate (k3) at (-1.3,0.8);
      \coordinate (k4) at (1.3,0.8);
      \coordinate (k5) at (0,-1.3);

      \fill[blue!15]
        (k1) -- (k2) -- (k3) -- cycle;
      \fill[blue!15]
        (k1) -- (k2) -- (k4) -- cycle;

      \foreach \i/\j in {1/2,1/3,2/3,1/4,2/4,3/5,4/5}
        \draw (k\i) -- (k\j);

      \foreach \i in {1,...,5}
        \node[vertex] at (k\i) {$\i$};

      \node at (0,-1.95)
        {(b) Clique complex $\operatorname{Cl}(G)$};
    \end{scope}

  \end{tikzpicture}
  \caption{A graph with five vertices and seven edges, and its
  clique complex. The two shaded triangles represent the
  two-dimensional simplices $\{1,2,3\}$ and $\{1,2,4\}$.
  The lower four-cycle is not filled.}
  \label{fig:clique-complex-example}
\end{figure}

For a concrete example, consider the graph $G$ on $V=[5]$ with
edge set
\[
  E=\{12,13,23,14,24,35,45\},
\]
where $ij$ abbreviates the unordered pair $\{i,j\}$.
Its only three-vertex cliques are $\{1,2,3\}$ and $\{1,2,4\}$,
and it has no larger clique. Consequently,
\[
  \operatorname{Cl}(G)
  =
  \{\varnothing\}
  \cup \bigl\{\{i\}:i\in[5]\bigr\}
  \cup E
  \cup \bigl\{\{1,2,3\},\{1,2,4\}\bigr\}.
\]
Thus $\operatorname{Cl}(G)$ has five vertices, seven edges and
two two-dimensional simplices. Geometrically, it consists of
two filled triangles sharing the edge $\{1,2\}$, together with
the edges $\{3,5\}$ and $\{4,5\}$; see
\cref{fig:clique-complex-example}.
The four-cycle $1,3,5,4,1$ remains unfilled.

Readers interested in clique complexes and graph cohomology
may consult the exposition in \cite{Lim2020}.

\begin{definition}
\label{def:induced-subgraph}
Let $G=(V,E)$ be a graph and let $S\subseteq V$. The
\emph{subgraph induced by $S$}, denoted by $G[S]$, is the graph
with vertex set $S$ and edge set
\[
  E(G[S])=\bigl\{\{u,v\}\in E:u,v\in S\bigr\}.
\]
Thus two vertices of $S$ are adjacent in $G[S]$ if and only if
they are adjacent in $G$. Equivalently, $G[S]$ is obtained by
deleting all vertices outside $S$ and all edges incident to them,
while retaining every edge whose endpoints both lie in $S$.
\end{definition}

\begin{figure}[htbp]
  \centering
  \begin{tikzpicture}[
    line width=0.7pt,
    line cap=round,
    line join=round,
    font=\small,
    vertex/.style={
      circle,
      draw=blue!60!black,
      fill=blue!12,
      minimum size=5mm,
      inner sep=0pt
    }
  ]
    \begin{scope}
      \coordinate (g1) at (0,0);
      \coordinate (g2) at (0,1.6);
      \coordinate (g3) at (-1.3,0.8);
      \coordinate (g4) at (1.3,0.8);
      \coordinate (g5) at (0,-1.3);

      \foreach \i/\j in {1/2,2/3,2/4}
        \draw[gray,dashed] (g\i) -- (g\j);
      \foreach \i/\j in {1/3,1/4,3/5,4/5}
        \draw[blue!60!black] (g\i) -- (g\j);

      \foreach \i in {1,3,4,5}
        \node[vertex] at (g\i) {$\i$};
      \node[vertex,draw=gray,fill=white,text=gray]
        at (g2) {$2$};

      \node at (0,-1.95) {(a) $G$, with $S$ highlighted};
    \end{scope}

    \node at (2.6,0.2) {$\longmapsto$};

    \begin{scope}[xshift=5.2cm]
      \coordinate (s1) at (0,0);
      \coordinate (s3) at (-1.3,0.8);
      \coordinate (s4) at (1.3,0.8);
      \coordinate (s5) at (0,-1.3);

      \foreach \i/\j in {1/3,1/4,3/5,4/5}
        \draw[blue!60!black] (s\i) -- (s\j);
      \foreach \i in {1,3,4,5}
        \node[vertex] at (s\i) {$\i$};

      \node at (0,-1.95) {(b) Induced subgraph $G[S]$};
    \end{scope}
  \end{tikzpicture}
  \caption{The subgraph induced by $S=\{1,3,4,5\}$.
  The highlighted vertices and edges are retained.
  The gray vertex and dashed edges belong to $G$ but are
  deleted when passing to $G[S]$.}
  \label{fig:induced-subgraph-example}
\end{figure}
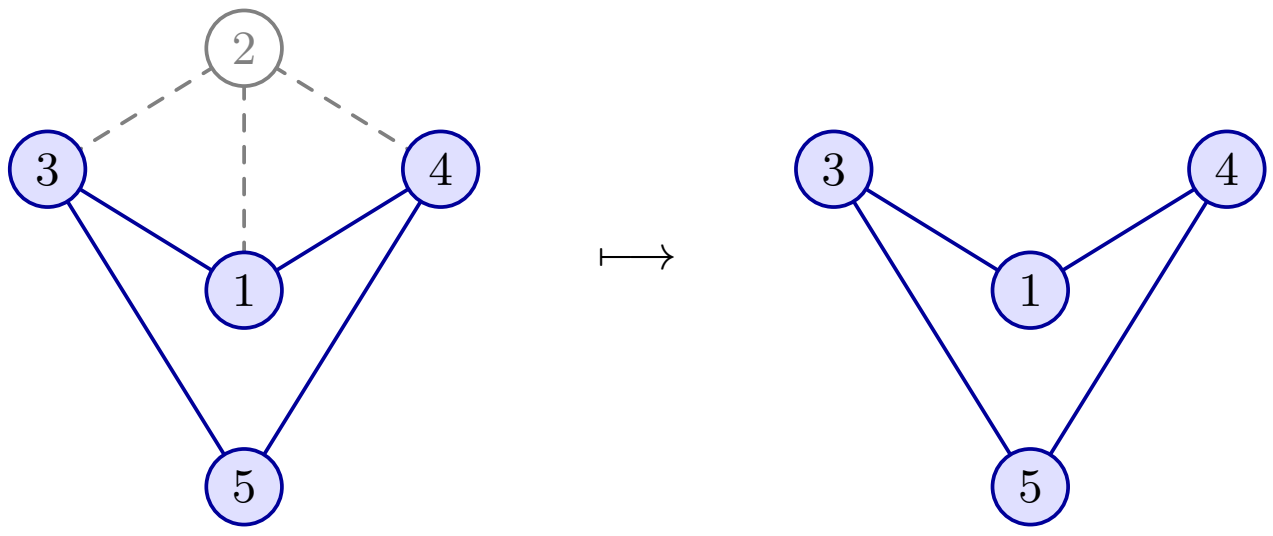

For example, let $G$ have vertex set $[5]$ and edge set
\[
  E=\{12,13,23,14,24,35,45\},
\]
where $ij$ abbreviates the unordered pair $\{i,j\}$.
For $S=\{1,3,4,5\}$, the induced subgraph has
\[
  V(G[S])=\{1,3,4,5\},
  \qquad
  E(G[S])=\{13,14,35,45\}.
\]
Deleting vertex $2$ removes precisely the edges $12$, $23$, and
$24$. The remaining edges form the cycle $1,3,5,4,1$;
see \cref{fig:induced-subgraph-example}.

\begin{definition}
\label{def:cycles-trees}
A \emph{path} of length $\ell\ge0$ is a subgraph with distinct
vertices $v_0,\ldots,v_\ell$ and precisely the edges
$\{v_{i-1},v_i\}$ for $1\le i\le\ell$. 
A \emph{cycle} of length $m\ge3$ is a subgraph with distinct
vertices $v_1,\ldots,v_m$ and precisely the edges
$\{v_1,v_2\},\ldots,\{v_{m-1},v_m\},\{v_m,v_1\}$.
A \emph{tree} is a connected graph with no cycle.
\end{definition}

\begin{definition}
\label{def:complement-graph}
Let $G=(V,E)$ be a finite simple graph. Its \emph{complement},
denoted by $\overline G$, is the graph with the same vertex set
$V$ and edge set
\[
  E(\overline G)
  =\bigl\{\{u,v\}:u,v\in V,\ u\ne v,\ \{u,v\}\notin E\bigr\}.
\]
Thus two distinct vertices are adjacent in $\overline G$ if
and only if they are not adjacent in $G$.
For $S\subseteq V$, the notation $\overline{G[S]}$ means the
complement of the induced subgraph $G[S]$, taken on the vertex
set $S$. In particular,
\[
  \overline{G[S]}=\overline G[S].
\]
\end{definition}

We next collect the spectral facts used in \cref{sg:section},
together with the form of Ramsey's theorem needed for its coverage
argument. In the spectral statements below, $J$ is a finite simple
graph with at least one vertex, and $M_J$ is its adjacency matrix.

\subsubsection{The adjacency spectrum}

The \emph{spectrum} of $J$ is the multiset of eigenvalues of $M_J$,
counted with their algebraic multiplicities. Since $M_J$ is real
and symmetric, these eigenvalues are real and may be ordered as
\[
  \lambda_1(M_J)\ge\cdots\ge\lambda_N(M_J),
  \qquad N=|V(J)|.
\]
We write $\lambda_{\min}(M_J)=\lambda_N(M_J)$ and define the
\emph{spectral radius} by
\[
  \rho(J):=\max_{1\le i\le N}|\lambda_i(M_J)|.
\]
For an adjacency matrix, $\rho(J)=\lambda_1(M_J)$, as follows
from the Perron--Frobenius theorem recalled below. Relabeling
vertices conjugates $M_J$ by a permutation matrix, so the
spectrum is independent of the vertex ordering. If $J$ is
$d$-regular, then $M_J\mathbf 1=d\mathbf 1$ and $\rho(J)=d$.

For a real symmetric matrix $A$ and a nonzero real vector $x$,
the \emph{Rayleigh quotient} is $x^TAx/(x^Tx)$. Orthogonal
diagonalization gives
\begin{equation}
  \lambda_{\min}(A)
  =\min_{x\ne0}\frac{x^TAx}{x^Tx},
  \qquad
  \lambda_{\max}(A)
  =\max_{x\ne0}\frac{x^TAx}{x^Tx}.
  \label{eq:prelim-rayleigh}
\end{equation}
A real symmetric matrix is \emph{positive semidefinite}, written
$A\succeq0$, if $x^TAx\ge0$ for every $x$. Equivalently, all
its eigenvalues are nonnegative. In particular,
\[
  \lambda_{\min}(M_J)\ge-3
  \quad\Longleftrightarrow\quad M_J+3I\succeq0.
\]
If $m_{-3}(J)$ denotes the multiplicity of $-3$ as an eigenvalue
of $M_J$, with value zero when $-3$ is absent, then
\[
  \rank(M_J+3I)=N-m_{-3}(J).
\]
Thus the normalized rank in \cref{sg:theorem} records the
proportion of adjacency eigenvalues different from $-3$.
We also use the Gram representation of a positive semidefinite
matrix: $A\succeq0$ if and only if there are real vectors $z_v$
such that $A_{uv}=\langle z_u,z_v\rangle$. This follows from a Cholesky 
factorization $A=ZZ^T$, with the rows of $Z$ as the vectors.

\subsubsection{Interlacing and forbidden induced subgraphs}

\begin{theorem}[{\cite[Corollary~2.5.2]{BH2012}}]
\label{thm:prelim-interlacing}
Let $A$ be a real symmetric $N\times N$ matrix, and let $B$ be
an $m\times m$ principal submatrix. With eigenvalues listed in
nonincreasing order,
\[
  \lambda_i(A)\ge\lambda_i(B)\ge\lambda_{N-m+i}(A)
  \qquad(1\le i\le m).
\]
\end{theorem}

For a nonempty vertex set $S\subseteq V(J)$, the adjacency
matrix of $J[S]$ is the principal submatrix $M_J[S]$. Hence
\begin{equation}
  \lambda_{\min}(M_J)
  \le\lambda_{\min}(M_{J[S]}).
  \label{eq:prelim-induced-interlacing}
\end{equation}
Consequently, if $\lambda_{\min}(M_J)\ge-\kappa$, then $J$
has no induced subgraph isomorphic to a graph $F$ with
$\lambda_{\min}(M_F)<-\kappa$. This is the exclusion principle
used in \cref{subsec:hoffman}.

For example, the \emph{star} $K_{1,m}$ consists of one center
joined to $m$ leaves, with no edges between the leaves. For
$m\ge1$, its eigenvalues are $\sqrt m$, $-\sqrt m$, and $0$
with multiplicity $m-1$. Indeed, its adjacency matrix is
\[
  \begin{pmatrix}0&\mathbf 1^T\\ \mathbf 1&0\end{pmatrix}.
\]
The subspace of vectors constant on the leaves gives the two
nonzero eigenvalues; vectors supported on the leaves with
coordinate sum zero form the kernel. Therefore
\begin{equation}
  \lambda_{\min}(M_J)\ge-3
  \quad\Longrightarrow\quad
  J\text{ has no induced }K_{1,10}.
  \label{eq:prelim-no-large-star}
\end{equation}
The same principle excludes $\widetilde K_{24}$ once a
Rayleigh quotient below $-3$ has been exhibited; this calculation
is given in \cref{subsec:hoffman}.

\subsubsection{Perron--Frobenius theory and subgraphs}

A \emph{subgraph} of $J$ may be obtained by deleting vertices
and edges. An induced subgraph retains every edge between its
remaining vertices. We recall the following result:

\begin{theorem}[{\cite[Proposition~3.1.1]{BH2012}}]
\label{thm:prelim-perron}
The matrix $M_J$ has a nonzero eigenvector with nonnegative
coordinates for the eigenvalue $\rho(J)=\lambda_1(M_J)$.
If $J$ is connected, this eigenvalue is simple and the
eigenvector can be chosen strictly positive.

If $H$ is a nonempty subgraph of $J$, then $\rho(H)\le\rho(J)$.
If, in addition, $J$ is connected and $H$ is a proper subgraph,
then $\rho(H)<\rho(J)$.
\end{theorem}

The subgraph comparison allows edge deletions, so $H$ need
not be induced. In particular, if $\rho(J)\le2$ and a
connected component of $J$ contains a connected subgraph
with spectral radius $2$, that component must equal the subgraph.
This observation will be used in \cref{sg:radius}.

\subsubsection{Smith's classification}

Write $P_n$ for the path on $n$ vertices and $C_n$ for the
cycle on $n\ge3$ vertices. A three-armed tree of arm lengths $(a,b,c)$
is formed by joining three otherwise disjoint paths of
$a$, $b$, and $c$ edges at a common endpoint. It has
$1+a+b+c$ vertices.

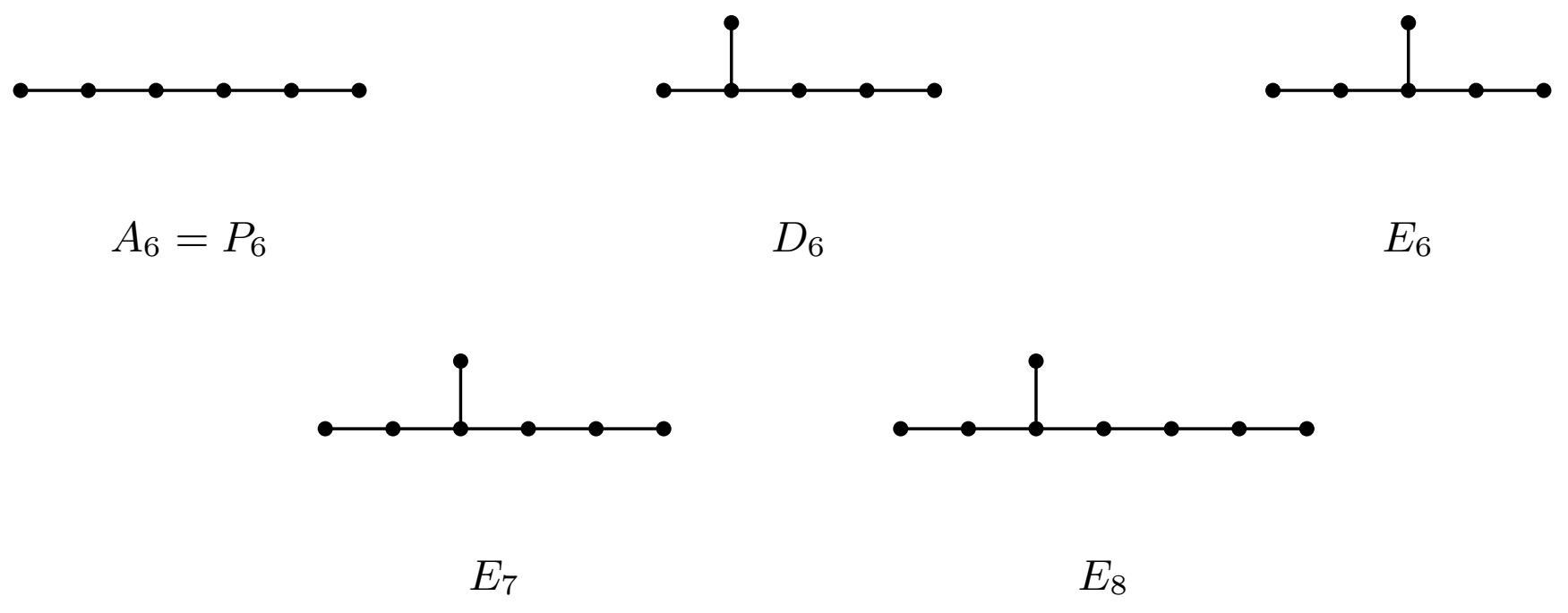
\begin{figure}[tbp]
  \centering
  \begin{tikzpicture}[
    x=0.5cm,y=0.5cm,line width=0.65pt,
    line cap=round,line join=round,font=\small]
    \begin{scope}
      \draw (-2.5,0)--(2.5,0);
      \foreach \x in {-2.5,-1.5,-0.5,0.5,1.5,2.5}
        \fill (\x,0) circle[radius=1.6pt];
      \node at (0,-2.2) {$A_6=P_6$};
    \end{scope}
    \begin{scope}[shift={(9,0)}]
      \draw (-2,0)--(2,0) (-1,0)--(-1,1);
      \foreach \x in {-2,-1,0,1,2}
        \fill (\x,0) circle[radius=1.6pt];
      \fill (-1,1) circle[radius=1.6pt];
      \node at (0,-2.2) {$D_6$};
    \end{scope}
    \begin{scope}[shift={(18,0)}]
      \draw (-2,0)--(2,0) (0,0)--(0,1);
      \foreach \x in {-2,-1,0,1,2}
        \fill (\x,0) circle[radius=1.6pt];
      \fill (0,1) circle[radius=1.6pt];
      \node at (0,-2.2) {$E_6$};
    \end{scope}
    \begin{scope}[shift={(4.5,-5)}]
      \draw (-2.5,0)--(2.5,0) (-0.5,0)--(-0.5,1);
      \foreach \x in {-2.5,-1.5,-0.5,0.5,1.5,2.5}
        \fill (\x,0) circle[radius=1.6pt];
      \fill (-0.5,1) circle[radius=1.6pt];
      \node at (0,-2.2) {$E_7$};
    \end{scope}
    \begin{scope}[shift={(13.5,-5)}]
      \draw (-3,0)--(3,0) (-1,0)--(-1,1);
      \foreach \x in {-3,-2,-1,0,1,2,3}
        \fill (\x,0) circle[radius=1.6pt];
      \fill (-1,1) circle[radius=1.6pt];
      \node at (0,-2.2) {$E_8$};
    \end{scope}
  \end{tikzpicture}
  \caption{The connected graphs with spectral radius less than
  $2$: representatives of the families $A_n$ and $D_n$, and
  the three exceptional trees.}
  \label{fig:smith-finite}
\end{figure}

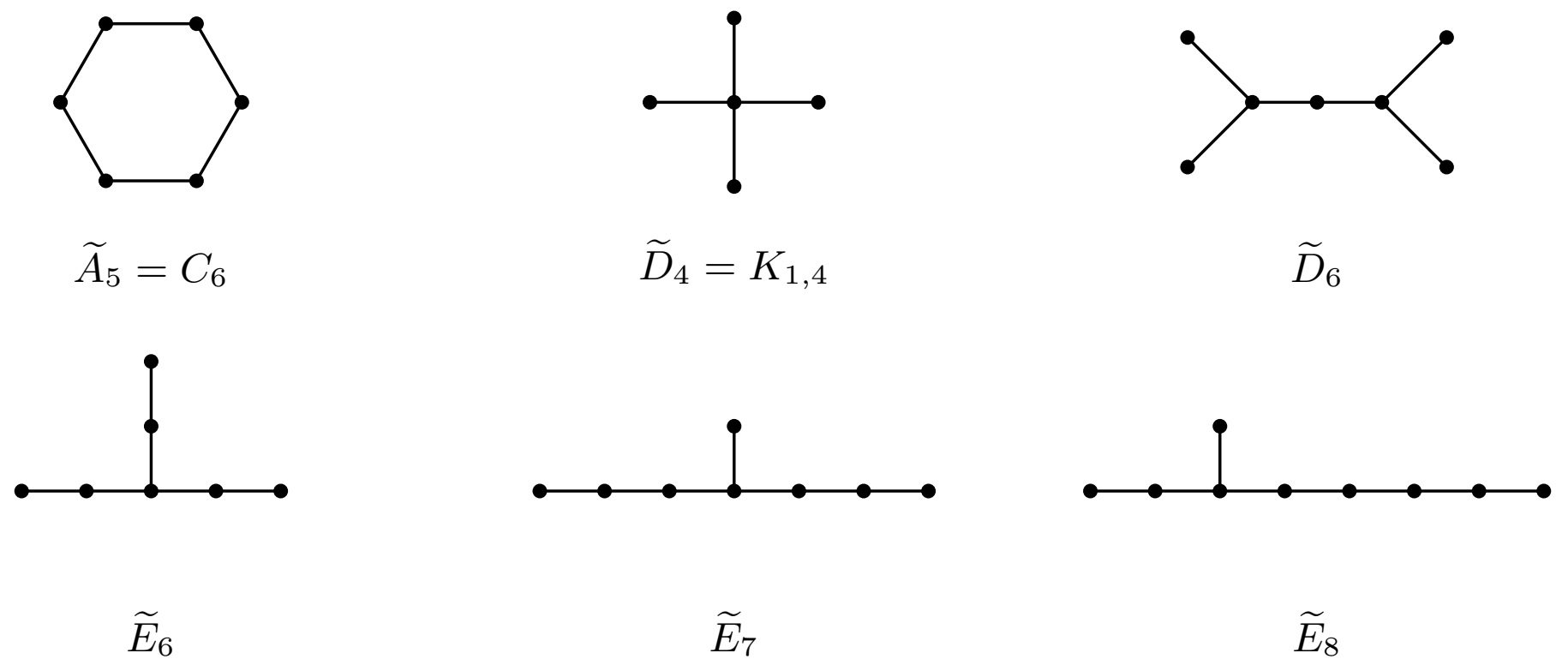
\begin{figure}[tbp]
  \centering
  \begin{tikzpicture}[
    x=0.5cm,y=0.5cm,line width=0.65pt,
    line cap=round,line join=round,font=\small]
    \begin{scope}
      \foreach \i in {0,...,5}
        \coordinate (c\i) at ({60*\i}:1.4);
      \draw (c0)--(c1)--(c2)--(c3)--(c4)--(c5)--cycle;
      \foreach \i in {0,...,5}
        \fill (c\i) circle[radius=1.6pt];
      \node at (0,-2.5) {$\widetilde A_5=C_6$};
    \end{scope}
    \begin{scope}[shift={(9,0)}]
      \draw (-1.3,0)--(1.3,0) (0,-1.3)--(0,1.3);
      \foreach \x/\y in {0/0,-1.3/0,1.3/0,0/-1.3,0/1.3}
        \fill (\x,\y) circle[radius=1.6pt];
      \node at (0,-2.5) {$\widetilde D_4=K_{1,4}$};
    \end{scope}
    \begin{scope}[shift={(18,0)}]
      \draw (-1,0)--(1,0)
        (-2,1)--(-1,0)--(-2,-1)
        (2,1)--(1,0)--(2,-1);
      \foreach \x/\y in {-2/1,-2/-1,-1/0,0/0,1/0,2/1,2/-1}
        \fill (\x,\y) circle[radius=1.6pt];
      \node at (0,-2.5) {$\widetilde D_6$};
    \end{scope}
    \begin{scope}[shift={(0,-6)}]
      \draw (-2,0)--(2,0) (0,0)--(0,2);
      \foreach \x in {-2,-1,0,1,2}
        \fill (\x,0) circle[radius=1.6pt];
      \foreach \y in {1,2}
        \fill (0,\y) circle[radius=1.6pt];
      \node at (0,-2.2) {$\widetilde E_6$};
    \end{scope}
    \begin{scope}[shift={(9,-6)}]
      \draw (-3,0)--(3,0) (0,0)--(0,1);
      \foreach \x in {-3,-2,-1,0,1,2,3}
        \fill (\x,0) circle[radius=1.6pt];
      \fill (0,1) circle[radius=1.6pt];
      \node at (0,-2.2) {$\widetilde E_7$};
    \end{scope}
    \begin{scope}[shift={(18,-6)}]
      \draw (-3.5,0)--(3.5,0) (-1.5,0)--(-1.5,1);
      \foreach \x in {-3.5,-2.5,-1.5,-0.5,0.5,1.5,2.5,3.5}
        \fill (\x,0) circle[radius=1.6pt];
      \fill (-1.5,1) circle[radius=1.6pt];
      \node at (0,-2.2) {$\widetilde E_8$};
    \end{scope}
  \end{tikzpicture}
  \caption{The connected graphs with spectral radius $2$:
  a cycle, the exceptional member $\widetilde D_4$ and a
  representative of the remaining $\widetilde D_n$ family,
  and the three extended exceptional trees.}
  \label{fig:smith-affine}
\end{figure}

The graphs in Smith's classification have the following
explicit descriptions. Put $A_n=P_n$ for $n\ge1$.
For $n\ge4$, let $D_n$ be the three-armed tree with arm
lengths $(1,1,n-3)$. The exceptional trees have arm lengths
\[
  E_6:(1,2,2),\qquad E_7:(1,2,3),\qquad E_8:(1,2,4).
\]
For the extended families, put $\widetilde A_{n-1}=C_n$
for $n\ge3$ and $\widetilde D_4=K_{1,4}$. For $n\ge5$,
construct $\widetilde D_n$ from a path of $n-4$ edges
by attaching two new leaves to each endpoint. Finally,
the arm lengths of the three extended exceptional trees are
\[
  \widetilde E_6:(2,2,2),\qquad
  \widetilde E_7:(1,3,3),\qquad
  \widetilde E_8:(1,2,5).
\]
Thus $A_n,D_n,E_n$ have $n$ vertices, whereas
$\widetilde A_n,\widetilde D_n,\widetilde E_n$ have $n+1$.
Examples of the families and all exceptional trees are
drawn in \cref{fig:smith-finite,fig:smith-affine}.

\begin{theorem}[{\cite[Theorem~3.1.3]{BH2012}}]
\label{thm:smith-classification}
Let $J$ be connected. Then $\rho(J)<2$ if and only if $J$
is isomorphic to one of
\[
  A_n\ (n\ge1),\qquad D_n\ (n\ge4),\qquad E_6,E_7,E_8.
\]
Moreover, $\rho(J)=2$ if and only if $J$ is isomorphic to one of
\[
  C_n\ (n\ge3),\qquad
  \widetilde D_n\ (n\ge4),\qquad
  \widetilde E_6,\widetilde E_7,\widetilde E_8.
\]
\end{theorem}

 Each tree
in the list has a path containing all but at most two vertices:
in a three-armed tree, use the two longest arms; in
$\widetilde D_n$, join a leaf at one end to a leaf at the other.
For $\widetilde D_4$, any path through the center and two leaves
has the required property. Paths in trees are induced, so this
 gives the induced paths used in \cref{sg:radius}.

\subsubsection{The Laplacian and connected components}

Let $D_J=\operatorname{diag}(\deg_J(v):v\in V(J))$.
The \emph{combinatorial Laplacian} of $J$ is $L_J=D_J-M_J$.
For every $x\in\R^{V(J)}$,
\begin{equation}
  x^TL_Jx
  =\sum_{\{u,v\}\in E(J)}(x_u-x_v)^2.
  \label{eq:prelim-laplacian}
\end{equation}
It follows that $L_J\succeq0$ and that its kernel consists
exactly of the vectors constant on each connected component.
If these components have vertex sets $W_1,\ldots,W_c$, then
\[
  \frac{\mathbf 1_{W_1}}{\sqrt{|W_1|}},\ldots,
  \frac{\mathbf 1_{W_c}}{\sqrt{|W_c|}}
\]
form an orthonormal basis of $\ker L_J$; here $\mathbf 1_W$
denotes the indicator vector of $W$. In particular,
$\dim\ker L_J=c$. These vectors are used to construct the
edge form in \cref{sg:kernel}.

\subsubsection{Ramsey's theorem}

A set of vertices is \emph{independent} if no two of them
are adjacent. For positive integers $s,t$, the \emph{Ramsey
number} $\mathrm R(s,t)$ is the least integer $N$ such that
every graph on at least $N$ vertices contains either a clique
of order $s$ or an independent set of order $t$.

\begin{theorem}[{\cite[Section~9.1]{Diestel2017}}]
\label{thm:prelim-ramsey}
The number $\mathrm R(s,t)$ is finite for all positive
integers $s,t$.
\end{theorem}

Write $N_J(v)=\{u\in V(J):uv\in E(J)\}$ for the neighborhood
of $v$. If $J$ has no induced $K_{1,t}$, then $N_J(v)$
contains no independent set of order $t$. Consequently,
every subset of $N_J(v)$ with at least $\mathrm R(s,t)$
vertices contains an $s$-clique. In \cref{sg:coverage}, we
apply this observation with $t=10$ after excluding
$K_{1,10}$ by interlacing.

For further background, we refer to \cite{BH2012} for graph
spectra and to \cite{Diestel2017} for Ramsey theory.

\subsection{Simplicial chains and cochains}
\label{subsec:chains-cochains}

We describe how to assign coefficients to oriented simplices
and how their boundaries and coboundaries are computed.
The resulting chain--cochain pairing and discrete Stokes formula
provide the algebraic framework for treating edge values as
cochains.

Let $K=\operatorname{Cl}(G)$, and work throughout with real
coefficients. A \emph{$q$-chain} on $K$ is a formal linear
combination of its $q$-dimensional simplices. To write these
combinations consistently, we list the vertices of each simplex
in increasing order and use the symbols $[i_0,\ldots,i_q]$ to denote the clique, 
with $i_0<\cdots<i_q$, as basis elements. Thus a $q$-chain has
the form
\[
  c=
  \sum_{\substack{i_0<\cdots<i_q\\
                   \{i_0,\ldots,i_q\}\in K}}
  c_{i_0\cdots i_q}[i_0,\ldots,i_q],
  \qquad c_{i_0\cdots i_q}\in\mathbb R.
\]
We denote the resulting vector space by
$C_q(K;\mathbb R)$. In particular, $0$-chains are weighted sums
of vertices, $1$-chains are weighted sums of edges, and
$2$-chains are weighted sums of triangles.

The ordering of the vertices records an \emph{orientation}.
For any permutation $\pi$ of $\{0,\ldots,q\}$, we impose the
sign convention
\[
  [i_{\pi(0)},\ldots,i_{\pi(q)}]
  =
  \operatorname{sgn}(\pi)[i_0,\ldots,i_q].
\]
Consequently, exchanging two vertices reverses the sign.
For example,
\[
  [j,i]=-[i,j],
  \qquad
  [j,k,i]=[i,j,k].
\]
The increasing vertex order therefore supplies a reference
orientation for each simplex. The underlying graph remains
undirected; orientation determines the signs used in the
chain spaces.

The \emph{boundary operator}
$\partial_q:C_q(K;\mathbb R)\to C_{q-1}(K;\mathbb R)$
is defined, for $q\ge1$, by
\begin{equation}
  \partial_q[i_0,\ldots,i_q]
  =
  \sum_{r=0}^{q}(-1)^r
  [i_0,\ldots,\widehat{i_r},\ldots,i_q],
  \label{eq:simplicial-boundary}
\end{equation}
and extended linearly to all chains.
Here the hat means that the indicated vertex is omitted.
We set $C_{-1}(K;\mathbb R)=0$ and $\partial_0=0$.
For an edge and a triangle, the formula becomes
\begin{align*}
  \partial_1[i,j]
    &= [j]-[i],\\
  \partial_2[i,j,k]
    &= [j,k]-[i,k]+[i,j]\\
    &= [i,j]+[j,k]+[k,i].
\end{align*}
Thus the boundary of an oriented edge is its terminal vertex
minus its initial vertex, while the boundary of an oriented
triangle follows its three edges around the triangle.

Applying the boundary operator twice gives zero:
\[
  \partial_{q-1}\partial_q=0.
\]
Indeed, each face obtained by deleting two vertices occurs
twice with opposite signs. A chain $c$ satisfying
$\partial_q c=0$ is called a \emph{$q$-cycle}, and a chain
of the form $c=\partial_{q+1}b$ is called a
\emph{$q$-boundary}. We write
\[
  Z_q(K;\mathbb R)=\ker\partial_q,
  \qquad
  B_q(K;\mathbb R)=\operatorname{im}\partial_{q+1}.
\]
The identity $\partial^2=0$ implies
$B_q(K;\mathbb R)\subseteq Z_q(K;\mathbb R)$:
every boundary is a cycle.

A \emph{$q$-cochain} assigns a real number to each oriented
$q$-simplex, with the same sign change under reversal of
orientation. Extending these assignments linearly to chains
identifies the cochain space with the dual vector space
\[
  C^q(K;\mathbb R)
  =
  \operatorname{Hom}\bigl(C_q(K;\mathbb R),\mathbb R\bigr).
\]
For $a\in C^q(K;\mathbb R)$ and
$c=\sum_\sigma c_\sigma\sigma\in C_q(K;\mathbb R)$,
their \emph{pairing} is
\[
  \langle a,c\rangle
  :=a(c)
  =\sum_\sigma c_\sigma a(\sigma),
\]
where the sum runs over the reference-oriented $q$-simplices.
A $0$-cochain is simply a function on the vertices.
A $1$-cochain is an assignment of edge values
$a_{ij}:=a([i,j])$ satisfying $a_{ji}=-a_{ij}$.

The \emph{coboundary operator}
$\delta_q:C^q(K;\mathbb R)\to C^{q+1}(K;\mathbb R)$
is defined by
\begin{equation}
  \langle\delta_q a,c\rangle
  =
  \langle a,\partial_{q+1}c\rangle,
  \qquad c\in C_{q+1}(K;\mathbb R).
  \label{eq:discrete-stokes}
\end{equation}
This identity is the \emph{discrete Stokes formula}:
evaluating a coboundary on a chain amounts to evaluating
the original cochain on its boundary. On a single simplex,
it reads
\[
  (\delta_q a)([i_0,\ldots,i_{q+1}])
  =
  \sum_{r=0}^{q+1}(-1)^r
  a([i_0,\ldots,\widehat{i_r},\ldots,i_{q+1}]).
\]
For a vertex function $f$ and an edge cochain $a$, respectively,
\begin{align*}
  (\delta_0 f)([i,j])
    &=f_j-f_i,\\
  (\delta_1 a)([i,j,k])
    &=a_{jk}-a_{ik}+a_{ij}\\
    &=a_{ij}+a_{jk}+a_{ki}.
\end{align*}
Thus $\delta_0$ takes differences of vertex values along
edges, while $\delta_1$ measures the oriented circulation
of edge values around each triangle.

The relation $\partial^2=0$ implies
$\delta_{q+1}\delta_q=0$. A cochain $a$ with $\delta_q a=0$
is called \emph{closed}, or a \emph{cocycle}; for $q\ge1$,
it is called \emph{exact} if $a=\delta_{q-1}b$ for some
$(q-1)$-cochain $b$. Every exact cochain is closed.
Exact cochains also vanish when paired with cycles:
for $b\in C^{q-1}(K;\mathbb R)$ and
$c\in Z_q(K;\mathbb R)$,
\[
  \langle\delta_{q-1}b,c\rangle
  =
  \langle b,\partial_q c\rangle
  =0.
\]
For an exact $1$-cochain, this is the familiar cancellation
of successive vertex differences around a closed walk. For further background on chain complexes, cochain complexes,
and homological algebra, we refer the reader to \cite{Weibel1994}.

\subsection{Exterior algebra, discrete products and Whitney forms}
\label{subsec:products-whitney}

We first recall the algebraic operations needed to multiply
forms and measure their sizes. We then describe their
simplexwise realization and the resulting products on cochains.
We conclude with classical helicity and its higher-dimensional
analogue, which motivate the local graph construction.

Readers concerned only with the discrete theory, rather than its
continuous motivation, may skip the discussion in
\cref{subsubsec:whitney-interpolation,subsubsec:cochain-products},
apart from \cref{prop:whitney-antisymmetrization}.
The normalized antisymmetrization formula in that proposition
may be taken as a purely discrete definition of the cochain
wedge product.

\subsubsection{Exterior algebra}
Let $U$ be an $N$-dimensional real Euclidean vector space.
A \emph{covector}, or algebraic one-form, is a linear map
$U\to\mathbb R$; the space of covectors is denoted by $U^*$.
More generally, $\Lambda^q U^*$ denotes the space of alternating
$q$-linear forms on $U$. 

Choose an orthonormal basis $e^1,\ldots,e^N$ of $U^*$.
The \emph{wedge product}, or exterior product, is bilinear
and associative, and satisfies
\[
  e^i\wedge e^j=-e^j\wedge e^i,
  \qquad
  e^i\wedge e^i=0.
\]
The products
\[
  e^{i_1}\wedge\cdots\wedge e^{i_q},
  \qquad i_1<\cdots<i_q,
\]
form a basis of $\Lambda^qU^*$.
Products can therefore be computed by concatenating basis
factors, reordering them with the appropriate signs, and
discarding terms containing a repeated factor.
For forms of degrees $p$ and $q$, respectively,
\[
  \alpha\wedge\gamma
  =
  (-1)^{pq}\gamma\wedge\alpha.
\]

The Euclidean structure on $U$ induces an inner product on
$\Lambda^qU^*$ by declaring the increasing exterior basis
above to be orthonormal. Thus
\[
  \alpha
  =
  \sum_{i_1<\cdots<i_q}
  \alpha_{i_1\cdots i_q}
  e^{i_1}\wedge\cdots\wedge e^{i_q}
  \quad\Longrightarrow\quad
  \|\alpha\|^2
  =
  \sum_{i_1<\cdots<i_q}\alpha_{i_1\cdots i_q}^{\,2}.
\]
This definition is independent of the orthonormal basis
chosen and fixes the norm convention used throughout the
Maclaurin argument.

A real skew-symmetric matrix $B=(b_{ij})$ determines an
algebraic two-form
\[
  \beta
  =
  \sum_{i<j}b_{ij}\,e^i\wedge e^j
  =
  \frac12\sum_{i,j}b_{ij}\,e^i\wedge e^j,
  \qquad
  \|\beta\|^2=\sum_{i<j}b_{ij}^2.
\]
We write $\beta^k$ for its $k$-fold wedge product and use
the normalized exterior power $\beta^k/k!$.
The factorial accounts for the different orderings of
commuting two-form factors. For example, if
\[
  \beta=x\,e^1\wedge e^2+y\,e^3\wedge e^4,
\]
then
\[
  \frac{\beta^2}{2!}
  =
  xy\,e^1\wedge e^2\wedge e^3\wedge e^4,
  \qquad
  \left\|\frac{\beta^2}{2!}\right\|^2=x^2y^2.
\]

For a more comprehensive treatment of tensor algebra and exterior
algebra, we refer the reader to \cite{Greub1978}.

\subsubsection{Differential forms on simplices}
A differential form allows the coefficients of an exterior
form to vary with the point. In affine coordinates
$x_1,\ldots,x_m$, a differential $q$-form has the expression
\[
  \alpha=\sum_{i_1<\cdots<i_q}
  f_{i_1\cdots i_q}(x)\,
  dx_{i_1}\wedge\cdots\wedge dx_{i_q}.
\]
Its \emph{exterior derivative} is defined by
\[
  d\alpha
  =
  \sum_{i_1<\cdots<i_q}
  df_{i_1\cdots i_q}
  \wedge dx_{i_1}\wedge\cdots\wedge dx_{i_q},
  \qquad
  df=\sum_{j=1}^{m}\frac{\partial f}{\partial x_j}\,dx_j.
\]
We reserve $d$ for this operator and $\delta$ for the
simplicial coboundary. The exterior derivative satisfies
\[
  d^2=0;
  \qquad
  d(\alpha\wedge\gamma)
  =
  d\alpha\wedge\gamma+(-1)^p\alpha\wedge d\gamma,
  \qquad \deg\alpha=p.
\]
Integration over an oriented simplex obeys Stokes' formula:
\[
  \int_\tau d\alpha=\int_{\partial\tau}\alpha.
\]
The integral over the boundary is the signed sum over its
oriented faces.

Every abstract $k$-simplex $\sigma=\{v_0,\ldots,v_k\}$ has a
standard Euclidean realization
\[
\Delta^k=\left\{(t_0,\ldots,t_k)\in\mathbb R^{k+1}:
 t_r\ge0,\ \sum_{r=0}^{k}t_r=1\right\},
\]
with $v_r$ represented by the coordinate unit vector $e_r$.
The unique coefficients in $x=\sum_{r=0}^{k}t_r e_r$ are its
\emph{barycentric coordinates}: as functions on this simplex,
$\lambda_{v_r}(x)=t_r$. Each is affine, equals one at its own
vertex and zero at the others, and
$\sum_{v\in\sigma}\lambda_v=1$, hence
$\sum_{v\in\sigma}d\lambda_v=0$.

A face is obtained by setting the coordinates of omitted vertices
to zero; the remaining coordinates restrict to its barycentric
coordinates. Euclidean models of simplices sharing a face are
glued by the affine map matching equally labelled vertices.
Differential forms on the complex are defined simplex by simplex,
with restrictions agreeing on common faces. This needs neither
a manifold structure nor an embedding of the whole complex in
a fixed Euclidean space. 
The ordering $[v_0,\ldots,v_k]$ orients the model by the basis
$(e_1-e_0,\ldots,e_k-e_0)$; exchanging two vertices reverses the
orientation. All simplex integrals and the Whitney interpolation
below are taken on these oriented Euclidean models, simplex by
simplex. Integration over a chain is the corresponding linear
combination of simplex integrals.

\makeatletter
\let\unmarkedseccntformat\@seccntformat
\renewcommand{\@seccntformat}[1]{*\unmarkedseccntformat{#1}}
\makeatother
\subsubsection{Whitney interpolation}
\label{subsubsec:whitney-interpolation}
Whitney forms convert cochain values into differential
forms on the simplices
\cite{Whitney1957,Hirani2003}.
For an oriented $q$-simplex
$\sigma=[i_0,\ldots,i_q]$, define
\begin{equation}
  W_\sigma
  =
  q!\sum_{r=0}^{q}(-1)^r\lambda_{i_r}\,
  d\lambda_{i_0}\wedge\cdots
  \wedge\widehat{d\lambda_{i_r}}\wedge\cdots
  \wedge d\lambda_{i_q}.
  \label{eq:prelim-whitney-form}
\end{equation}
On every simplex containing $\sigma$, the coordinates in
this formula are the barycentric coordinates of that
containing simplex. On simplices not containing $\sigma$,
the form is zero. The restrictions agree on common faces,
so $W_\sigma$ defines a simplexwise differential form on $K$.

For $a\in C^q(K;\mathbb R)$, its \emph{Whitney interpolation}
is
\[
  W_Ka=\sum_{\dim\sigma=q}a(\sigma)W_\sigma,
\]
where the sum uses the reference orientations.
For example, a vertex function is interpolated by
$W_Kf=\sum_i f_i\lambda_i$.
For an edge cochain $a$, the restriction to a simplex $\tau$
is
\begin{align*}
  (W_Ka)|_\tau
    &=
    \sum_{\substack{i<j\\\{i,j\}\subseteq\tau}}
    a_{ij}
    \bigl(\lambda_i\,d\lambda_j-\lambda_j\,d\lambda_i\bigr),\\
  (dW_Ka)|_\tau
    &=
    2\sum_{\substack{i<j\\\{i,j\}\subseteq\tau}}
    a_{ij}\,d\lambda_i\wedge d\lambda_j.
\end{align*}
Thus the edge values determine a differential one-form
throughout each simplex containing those edges.

Conversely, integration defines a map from compatible
simplexwise differential forms to cochains:
\[
  (R_K\alpha)(\sigma)=\int_\sigma\alpha.
\]
The normalization in \eqref{eq:prelim-whitney-form} ensures
that integration recovers the original cochain.
The fundamental compatibility identities are
\begin{equation}
  R_KW_K=\mathrm{id},
  \qquad
  dW_K=W_K\delta,
  \qquad
  R_Kd=\delta R_K.
  \label{eq:prelim-whitney-identities}
\end{equation}
The last identity is Stokes' formula expressed in terms
of the integration map. For further details on Whitney forms and their compatibility
with simplicial coboundaries and integration, we refer the
reader to \cite{Whitney1957,Hirani2003}.

There are two different uses of edge data in this paper.
An edge cochain $a$ gives a simplexwise one-form $W_Ka$.
It also gives, after setting $a_{ij}=0$ on nonedges, an
algebraic two-form
\[
  \beta_a=\sum_{i<j}a_{ij}\,e^i\wedge e^j
  \quad\text{on }\mathbb R^n.
\]
The latter identification will be used in the Maclaurin
argument. Its exterior-algebra norm is
$\|\beta_a\|^2=\sum_{\{i,j\}\in E,\ i<j}a_{ij}^2$.

\subsubsection{Products of cochains}
\label{subsubsec:cochain-products}
We use two products that combine a $p$-cochain and a
$q$-cochain into a $(p+q)$-cochain.
The \emph{cup product} is defined on each
increasingly ordered simplex by
\begin{equation}
  (a\smile b)
  ([i_0,\ldots,i_{p+q}])
  =
  a([i_0,\ldots,i_p])\,
  b([i_p,\ldots,i_{p+q}]),
  \qquad i_0<\cdots<i_{p+q},
  \label{eq:prelim-AW-product}
\end{equation}
and then extended according to the orientation convention
and linearity. For two edge cochains, for example,
\[
  (a\smile b)([i,j,k])
  =
  a_{ij}b_{jk},
  \qquad i<j<k.
\]
This product is associative, but is generally not
graded-commutative on cochains.

The \emph{Whitney product} of cochains is obtained by
interpolating, taking the wedge product, and integrating:
\begin{equation}
  a\wedge b
  :=
  R_K\bigl(W_Ka\wedge W_Kb\bigr).
  \label{eq:prelim-whitney-product}
\end{equation}
We use the same wedge notation for cochains and differential
forms, with its meaning determined by the arguments.
For $a\in C^p(K;\mathbb R)$ and $b\in C^q(K;\mathbb R)$,
\[
  a\wedge b=(-1)^{pq}b\wedge a.
\]
Both cochain products obey the Leibniz identity
\[
  \delta(a\diamond b)
  =
  (\delta a)\diamond b
  +
  (-1)^p a\diamond(\delta b),
  \qquad
  \diamond\in\{\smile,\wedge\}.
\]
The wedge product of cochains is generally not associative.
Powers of differential forms, such as $(dW_Ka)^s$, use the
associative exterior product.

The Whitney product also has a purely combinatorial expression:
it is the normalized antisymmetrization of the raw vertex
expression underlying the cup product.
\begin{proposition}
\label{prop:whitney-antisymmetrization}
Let $K$ be a simplicial complex, and let
$a\in C^p(K;\mathbb R)$ and $b\in C^q(K;\mathbb R)$.
For every oriented $(p+q)$-simplex $[v_0,\ldots,v_{p+q}]$,
the Whitney product satisfies
\begin{equation}
  \begin{aligned}
    (a\wedge b)([v_0,\ldots,v_{p+q}])
    &=
    \frac{1}{(p+q+1)!}
    \sum_{\pi\in S_{p+q+1}}\operatorname{sgn}(\pi)\\
    &\qquad{}\cdot
    a([v_{\pi(0)},\ldots,v_{\pi(p)}])\,
    b([v_{\pi(p)},\ldots,v_{\pi(p+q)}]).
  \end{aligned}
  \label{eq:prelim-wedge-antisymmetrization}
\end{equation}
Here $S_{p+q+1}$ denotes the group of permutations of
$\{0,\ldots,p+q\}$.

\end{proposition}

\begin{proof}
The combinatorial expression on the right-hand side is
the discrete wedge product in
\cite[Definition~7.2.1]{Hirani2003}.
Its equality with $R_K(W_Ka\wedge W_Kb)$ is proved in
\cite[Proposition~1]{SchubelBerwickEvansHirani2024}.
The proof is local to each simplex and therefore applies
to an arbitrary simplicial complex.
\end{proof}

\subsubsection{Classical helicity and its higher-dimensional analogue}

For a smooth one-form $\alpha$ on a closed oriented
Riemannian three-manifold $(M,g)$, the helicity of the exact
two-form $d\alpha$ is
\[
  \mathcal H(d\alpha)
  :=
  \int_M \alpha\wedge d\alpha.
\]
This value is independent of the choice of potential
$\alpha$. Although the definition requires only an orientation,
a Riemannian metric allows helicity to be compared with energy.
The classical Arnold energy inequality states that
\[
  |\mathcal H(d\alpha)|
  \le C_{M,g}\,\|d\alpha\|_{L^2(M)}^2,
\]
where $C_{M,g}>0$ depends on the manifold and its
Riemannian metric \cite{ArnoldKhesin2021}.
Thus nonzero helicity gives a positive lower bound for
the energy of the exact two-form.

The same differential-form construction extends to higher
odd dimensions \cite{Khesin1992,ArnoldKhesin2021}.
Let $s\ge2$, let $M$ be a closed oriented manifold of
dimension $2s-1$, and let $\beta=d\alpha$ be a smooth exact
two-form on $M$. Define
\begin{equation}
  \mathcal H_s(\beta)
  :=
  \int_M \alpha\wedge\beta^{s-1}
  =
  \int_M \alpha\wedge(d\alpha)^{s-1}.
  \label{eq:prelim-higher-helicity}
\end{equation}
Here powers denote repeated exterior products.
The integrand has degree $1+2(s-1)=2s-1$, so it can be
integrated over $M$. When $s=2$, this definition recovers
the classical three-dimensional helicity.

As in dimension three, $\mathcal H_s(\beta)$ is independent
of the choice of potential. Indeed, any other potential
has the form $\alpha+z$ with $dz=0$, and
\[
  z\wedge\beta^{s-1}
  =
  -d\bigl(z\wedge\alpha\wedge\beta^{s-2}\bigr).
\]
Since $M$ has no boundary, Stokes' theorem gives
$\int_M z\wedge\beta^{s-1}=0$.
Consequently, \eqref{eq:prelim-higher-helicity} defines
a functional of the exact two-form $\beta$ itself.

These formulas suggest a corresponding local construction
on a simplicial complex. Replacing $\alpha$ by the Whitney
interpolation $W_Ka$ of an edge cochain gives the
simplexwise differential form
\[
  W_Ka\wedge(dW_Ka)^{s-1},
\]
which can be integrated over each oriented
$(2s-1)$-simplex. In \cref{sec:helicity}, we apply this
construction to the simplices of a graph's clique complex
and show that, after an explicit normalization, each local
integral is exactly the Pfaffian of the corresponding
edge matrix.

For a broader account of helicity and its role in topological
hydrodynamics, we refer the reader to \cite{ArnoldKhesin2021}.

\makeatletter
\let\@seccntformat\unmarkedseccntformat
\makeatother
\subsection{Pfaffians and classical Maclaurin inequalities}
\label{subsec:pfaffians-maclaurin}

We recall the Pfaffian as a signed sum over perfect matchings
and explain how it records the coefficients of exterior powers.
We then state the classical Maclaurin inequalities that will
be used to compare these powers.

\subsubsection{Pfaffians}

A perfect matching of a set of $2m$ vertices is a collection
of $m$ disjoint pairs covering every vertex exactly once.
The Pfaffian assigns a signed weight to each such matching.

\begin{definition}
\label{def:pfaffian}
Let $B=(b_{ij})$ be a real skew-symmetric matrix of order $2m$,
and let $\mathcal M_{2m}$ denote the set of perfect matchings
of $[2m]$.
Write each matching uniquely as
\[
  M=\bigl\{\{i_1,j_1\},\ldots,\{i_m,j_m\}\bigr\},
  \qquad
  i_\ell<j_\ell,\quad i_1<\cdots<i_m.
\]
Let $\varepsilon(M)$ be the sign of the permutation
\[
  (1,2,\ldots,2m)\longmapsto
  (i_1,j_1,\ldots,i_m,j_m).
\]
The \emph{Pfaffian} of $B$ is
\begin{equation}
  \Pf(B)
  :=
  \sum_{M\in\mathcal M_{2m}}
  \varepsilon(M)\prod_{\ell=1}^{m}b_{i_\ell j_\ell}.
  \label{eq:prelim-pfaffian-definition}
\end{equation}
We set the Pfaffian of the empty matrix equal to one.
\end{definition}

For a matrix of order two, the definition gives
\[
  \Pf\begin{pmatrix}
    0&b_{12}\\
    -b_{12}&0
  \end{pmatrix}
  =b_{12}.
\]
For a matrix of order four, there are precisely three
perfect matchings, and
\begin{equation}
  \Pf(B)
  =
  b_{12}b_{34}-b_{13}b_{24}+b_{14}b_{23}.
  \label{eq:prelim-four-point-pfaffian}
\end{equation}
The three terms correspond to the matchings
$\{12,34\}$, $\{13,24\}$ and $\{14,23\}$, respectively.
If the entries of $B$ vanish on the nonedges of a graph,
only perfect matchings using its edges can contribute.
Different matching terms may nevertheless cancel.

The Pfaffian satisfies the classical identity
\begin{equation}
  \det B=\Pf(B)^2.
  \label{eq:prelim-pfaffian-determinant}
\end{equation}
In particular, an even-order skew-symmetric matrix is
singular if and only if its Pfaffian vanishes.
Unlike its determinant, the Pfaffian retains a sign:
applying the same permutation to the rows and columns
multiplies the Pfaffian by the sign of that permutation.

\subsubsection{Pfaffians and exterior powers}

Let $B=(b_{ij})$ now be a real skew-symmetric matrix of
order $N$, and associate to it the two-form
\[
  \beta=\sum_{1\le i<j\le N}b_{ij}\,e^i\wedge e^j,
\]
where $e^1,\ldots,e^N$ is an orthonormal basis of covectors.
For a subset $S=\{i_1<\cdots<i_{2s}\}\subseteq[N]$, write
$B[S]$ for the principal submatrix whose rows and columns
are indexed by $S$ in increasing order, and put
\[
  e^S=e^{i_1}\wedge\cdots\wedge e^{i_{2s}}.
\]
The Pfaffians of these principal submatrices are exactly
the coefficients of the normalized exterior power:
\begin{equation}
  \frac{\beta^s}{s!}
  =
  \sum_{\substack{S\subseteq[N]\\|S|=2s}}
  \Pf(B[S])\,e^S,
  \qquad
  1\le s\le\left\lfloor\frac N2\right\rfloor.
  \label{eq:prelim-pfaffian-exterior}
\end{equation}
Indeed, a nonzero term in the expansion of $\beta^s$ uses
$2s$ distinct indices and therefore determines a perfect
matching of those indices. Reordering its covectors gives
the sign in \eqref{eq:prelim-pfaffian-definition}.
The two-form factors commute, so each matching occurs
$s!$ times, accounting for the normalization.

Since the increasing exterior basis is orthonormal,
\eqref{eq:prelim-pfaffian-exterior} also gives
\begin{equation}
  \left\|\frac{\beta^s}{s!}\right\|^2
  =
  \sum_{\substack{S\subseteq[N]\\|S|=2s}}
  \Pf(B[S])^2.
  \label{eq:prelim-pfaffian-norm}
\end{equation}

\subsubsection{Elementary symmetric polynomials and Maclaurin inequalities}

\begin{definition}
\label{def:elementary-symmetric}
For real numbers $x_1,\ldots,x_r$ and $1\le s\le r$,
the \emph{$s$-th elementary symmetric polynomial} is
\[
  e_s(x_1,\ldots,x_r)
  :=
  \sum_{1\le i_1<\cdots<i_s\le r}
  x_{i_1}\cdots x_{i_s}.
\]
We also set $e_0=1$ and $e_s=0$ for $s>r$.
\end{definition}

Thus
\[
  e_1(x)=\sum_{i=1}^{r}x_i,
  \qquad
  e_2(x)=\sum_{i<j}x_ix_j,
  \qquad
  e_r(x)=\prod_{i=1}^{r}x_i.
\]
There are $\binom rs$ products in the definition of $e_s$.
Consequently, $e_s(x)/\binom rs$ is their average.
For nonnegative inputs, taking its $s$-th root allows
products of different degrees to be compared.

\begin{theorem}[Maclaurin inequalities {\cite{HLP1952}}]
\label{thm:classical-maclaurin}
Let $r\ge2$ and $x_1,\ldots,x_r\ge0$.
For every pair of integers $1\le j<k\le r$,
\begin{equation}
  \left(
    \frac{e_k(x_1,\ldots,x_r)}{\binom rk}
  \right)^{1/k}
  \le
  \left(
    \frac{e_j(x_1,\ldots,x_r)}{\binom rj}
  \right)^{1/j}.
  \label{eq:prelim-classical-maclaurin}
\end{equation}
Equivalently,
\[
  \frac{e_1(x)}{r}
  \ge
  \left(\frac{e_2(x)}{\binom r2}\right)^{1/2}
  \ge\cdots\ge
  e_r(x)^{1/r}.
\]
\end{theorem}

In \cref{lem:two-form-maclaurin}, we apply these inequalities
to the squared parameters of the $2\times2$ blocks in the
real skew-symmetric normal form.
Together with \eqref{eq:prelim-pfaffian-norm}, this relates
the classical symmetric-polynomial inequalities to sums
of squared principal Pfaffians.

\section{Local helicity and its Pfaffian realization}
\label{sec:helicity}

The purpose of this section is to express local Whitney helicity
entirely in terms of edge values. On a tetrahedron, the helicity
integral is one third of the Pfaffian of its edge matrix. More generally, on a $(2s-1)$-simplex, it is a dimension-dependent multiple of the corresponding Pfaffian. We determine this factor
and normalize the local helicities so that they are exactly
Pfaffians. Their squared values then define the graph energies
used in the Maclaurin inequality.

\subsection{Local helicity on cliques}
\label{subsec:local-helicity}

Let $K=\Cl(G)$, with the fixed ordering of its vertices.
We identify an edge cochain $a\in C^1(K;\R)$ with the
skew-symmetric matrix $(a_{ij})_{i,j\in V(G)}$ obtained by setting
$a_{ij}=a([i,j])$ on edges and $a_{ij}=0$ on nonedges and on the
diagonal. We abbreviate $C^1(K;\R)$ as $C^1(G)$. We write
$d_{\mathrm{dR}}$ for the exterior derivative in this section.

Every clique with an even number of vertices carries a local
helicity integral on its corresponding simplex of $K$. For
\[
  S=\{i_0<\cdots<i_{2s-1}\}\in\Kcal_{2s}(G),
\]
let $\Delta_S^{2s-1}$ be the corresponding affine simplex,
oriented by $[i_0,\ldots,i_{2s-1}]$, and let
$W_S(a)=(W_Ka)|_{\Delta_S^{2s-1}}$.
In barycentric coordinates $t_0,\ldots,t_{2s-1}$, this one-form is
\[
  W_S(a)
  =
  \sum_{0\le\alpha<\beta\le2s-1}
  a_{i_\alpha i_\beta}
  \bigl(t_\alpha\,dt_\beta-t_\beta\,dt_\alpha\bigr),
  \qquad
  \sum_{\alpha=0}^{2s-1}t_\alpha=1.
\]

\begin{definition}
\label{def:local-whitney-helicity}
For $s\ge1$, an edge form $a\in C^1(G)$, and
$S\in\Kcal_{2s}(G)$ with its reference orientation, define
\begin{equation}
  \mathcal H_{s,S}(a)
  :=
  \int_{\Delta_S^{2s-1}}
  W_S(a)\wedge\bigl(d_{\mathrm{dR}}W_S(a)\bigr)^{s-1}
  =
  \bigl(a\wedge(\delta a)^{s-1}\bigr)(S).
  \label{eq:local-whitney-helicity}
\end{equation}
Here $(\delta a)^0$ denotes the unit $0$-cochain.
Since $d_{\mathrm{dR}}W_S(a)$ has constant coefficients,
Whitney interpolation reproduces its exterior powers,
which justifies the second equality.
\end{definition}

For $s=2$, this recovers the tetrahedral helicity integral
in Kotiuga's Whitney-form discretization
\cite{Kotiuga1989}.
The construction requires only the simplex associated with
each clique and therefore applies to an arbitrary abstract
graph, without an embedding or a manifold structure.

\subsection{The Pfaffian formula and normalization}
\label{subsec:higher-local}

The integral in \eqref{eq:local-whitney-helicity} has an algebraic
expression with a universal coefficient. The following identity
identifies that coefficient in every dimension.

\begin{proposition}
\label{prop:higher-whitney}
For every $s\ge1$, every $S\in\Kcal_{2s}(G)$, and every
$a\in C^1(G)$,
\begin{equation}
  \mathcal H_{s,S}(a)
  =\frac{2^{s-1}(s-1)!}{(2s-1)!}\Pf\bigl(a[S]\bigr).
  \label{eq:local-whitney-pfaffian}
\end{equation}
Consequently, the normalized local helicity
\begin{equation}
  P_{s,S}(a)
  :=\frac{(2s-1)!}{2^{s-1}(s-1)!}\mathcal H_{s,S}(a)
  =\Pf\bigl(a[S]\bigr)
  \label{eq:normalized-local-helicity}
\end{equation}
is exactly the Pfaffian of the restricted edge matrix.
\end{proposition}

\begin{proof}
Relabel the vertices of \(S\) locally by \(0,\dots,2s-1\), and put
\[
\Omega
=
\sum_{\alpha<\beta}
a_{i_\alpha i_\beta}\,
dt_\alpha\wedge dt_\beta
\]
in the ambient barycentric space. Also let
\[
X
=
\sum_{\alpha=0}^{2s-1}
t_\alpha\frac{\partial}{\partial t_\alpha}
\]
be the Euler vector field.

First,
\[
d_{\mathrm{dR}}
\bigl(
t_\alpha\,dt_\beta-t_\beta\,dt_\alpha
\bigr)
=
2\,dt_\alpha\wedge dt_\beta,
\]
and hence
\[
d_{\mathrm{dR}}W_S(a)=2\Omega.
\]
Moreover,
\[
\iota_X(dt_\alpha\wedge dt_\beta)
=
t_\alpha\,dt_\beta-t_\beta\,dt_\alpha,
\]
so
\[
\iota_X\Omega=W_S(a).
\]

Since \(\Omega\) has even degree, the contraction Leibniz rule gives
\[
\iota_X(\Omega^s)
=
s(\iota_X\Omega)\wedge\Omega^{s-1}
=
s\,W_S(a)\wedge\Omega^{s-1}.
\]
Therefore
\[
W_S(a)\wedge
\bigl(d_{\mathrm{dR}}W_S(a)\bigr)^{s-1}
=
\frac{2^{s-1}}{s}\,
\iota_X(\Omega^s).
\]

We now expand \(\Omega^s\). A nonzero monomial must use every one of the
\(2s\) indices exactly once, hence corresponds to a perfect matching of
\(\{0,\dots,2s-1\}\). The signed sum over such perfect matchings is precisely
the Pfaffian. Since the \(s\) two-form factors commute with one another,
each matching appears \(s!\) times. Thus
\[
\Omega^s
=
s!\,\Pf\bigl(a[S]\bigr)\,
dt_0\wedge\cdots\wedge dt_{2s-1}.
\]
Consequently,
\[
W_S(a)\wedge
\bigl(d_{\mathrm{dR}}W_S(a)\bigr)^{s-1}
=
2^{s-1}(s-1)!
\Pf\bigl(a[S]\bigr)\,
\iota_X
\bigl(
dt_0\wedge\cdots\wedge dt_{2s-1}
\bigr).
\]

It remains to restrict the contracted ambient volume form to the affine simplex
\[
t_0+\cdots+t_{2s-1}=1.
\]
Expanding the contraction,
\[
\iota_X
\bigl(
dt_0\wedge\cdots\wedge dt_{2s-1}
\bigr)
=
\sum_{\alpha=0}^{2s-1}
(-1)^\alpha t_\alpha\,
dt_0\wedge\cdots\widehat{dt_\alpha}\cdots
\wedge dt_{2s-1}.
\]
On the affine hyperplane,
\[
dt_0=-\sum_{\alpha=1}^{2s-1}dt_\alpha.
\]
After substitution, every summand reduces with positive sign to
\[
t_\alpha\,dt_1\wedge\cdots\wedge dt_{2s-1}.
\]
Since \(\sum_\alpha t_\alpha=1\), the restriction is therefore
\[
dt_1\wedge\cdots\wedge dt_{2s-1}.
\]
Finally,
\[
\int_{\Delta_S^{2s-1}}
dt_1\wedge\cdots\wedge dt_{2s-1}
=
\frac1{(2s-1)!}.
\]
Combining the preceding identities proves the formula.
\end{proof}

For a tetrahedron $\sigma=[i,j,k,\ell]$ with $i<j<k<\ell$,
\cref{prop:higher-whitney} gives the explicit identity
\begin{equation}
  \mathcal H_{2,\sigma}(a)
  =\frac13\bigl(
    a_{ij}a_{k\ell}-a_{ik}a_{j\ell}+a_{i\ell}a_{jk}
  \bigr).
  \label{eq:four-point-local-helicity}
\end{equation}
In this case we abbreviate the normalized local helicity by
\begin{equation}
  P_\sigma(a):=P_{2,\sigma}(a)
  =3\mathcal H_{2,\sigma}(a)
  =\Pf\bigl(a[\sigma]\bigr).
  \label{eq:normalized-four-point-helicity}
\end{equation}

\IfFileExists{figures/pfaffian-matchings.pdf}{%
\begin{figure}[tb]
\centering
\includegraphics[width=.92\linewidth]{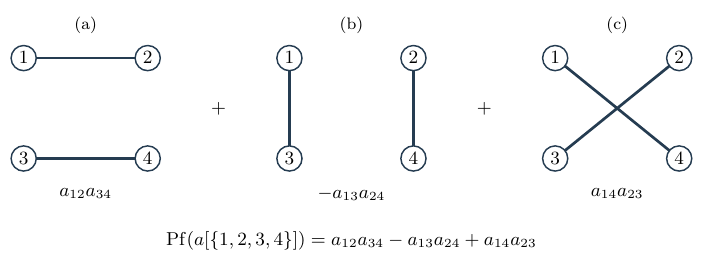}
\caption{The three perfect matchings of an ordered four-vertex
clique contribute the signs $+,-,+$ to its Pfaffian. The local
Whitney helicity integral is one third of their signed sum.}
\label{fig:pfaffian-matchings}
\end{figure}}{}

\subsection{From local helicity to graph energies}
\label{subsec:helicity-energies}

We aggregate the normalized local helicities by summing their
squares over cliques. This produces nonnegative quantities that
retain the magnitude of each local contribution even when signed
helicities cancel.

\begin{definition}
\label{def:pfaffian-energies}
For an edge form $a\in C^1(G)$ and an integer $s\ge1$, define
\begin{equation}
  T_{s,G}(a)
  :=\sum_{S\in\Kcal_{2s}(G)}P_{s,S}(a)^2
  =\sum_{S\in\Kcal_{2s}(G)}\Pf\bigl(a[S]\bigr)^2.
  \label{eq:energy-definition}
\end{equation}
An empty sum is zero. The first energy is the squared norm
of the edge form:
\[
  E_G(a):=T_{1,G}(a)
  =\sum_{\substack{\{i,j\}\in E\\i<j}}a_{ij}^2.
\]
\end{definition}

In particular, for $s=2$ the energy is
\[
  T_{2,G}(a)
  =\sum_{\substack{i<j<k<\ell\\
                 \{i,j,k,\ell\}\in\Kcal_4(G)}}
    \bigl(a_{ij}a_{k\ell}-a_{ik}a_{j\ell}+a_{i\ell}a_{jk}\bigr)^2
  =9\sum_{\sigma\in\Kcal_4(G)}\mathcal H_{2,\sigma}(a)^2.
\]
These energies are defined on every graph, including graphs whose
clique complexes have no nonzero three-cycles. They are also
unchanged by reordering the vertices, since each summand is a
squared Pfaffian.

The local helicities above do not, by themselves, determine a canonical scalar helicity for an arbitrary graph. On a closed
oriented manifold, local contributions can be integrated over
its fundamental cycle. An arbitrary clique complex need not
possess an analogous fundamental class. The fixed vertex order
provides reference orientations for individual simplices, but
does not ensure that their sum is a cycle. A global helicity
analogous to the manifold integral would therefore require
additional data specifying an appropriate cycle.

 In the next section,
we compare these energies at different orders and show that
their sharp Maclaurin normalization is determined by the
clique number of $G$.

\section{The helicity energy Maclaurin inequality}
\label{sec:maclaurin}

We now compare the energies $T_{s,G}$ defined in
\eqref{eq:energy-definition}. A closely related line of work
studies clique polynomials. Motzkin and
Straus~\cite{MotzkinStraus1965} characterized the clique number
through quadratic optimization over nonnegative vertex weights
summing to one.
Fisher and Ryan~\cite{FisherRyan1992} proved Maclaurin inequalities
for clique counts, and Nikiforov~\cite{Nikiforov2006} established
the corresponding inequalities for vertex-weighted clique
polynomials, normalized by $\binom{\omega(G)}{s}$.
Arag{\~a}o and Souza~\cite{AragaoSouza2024} obtained localised
versions whose weights depend on the largest clique containing
each contributing clique.

The energies considered here sum squared Pfaffians of an edge
form over $2s$-cliques, rather than monomials in vertex weights.
Their comparison has the form of a Maclaurin inequality with
normalization $\binom{\lfloor\omega(G)/2\rfloor}{s}$.
Our proof uses a related adjustment between nonadjacent vertices
to reduce the support to a clique, together with Pfaffian zero
propagation and the Euclidean two-form inequality.

\begin{maintheorem}
\label{thm:masked-maclaurin}
Let $G$ be a finite simple graph and put
$p=\lfloor\omega(G)/2\rfloor$. For every edge form
$a\in C^1(G)$ and all integers $1\le j<k\le p$,
\[
  \left(\frac{T_{k,G}(a)}{\binom{p}{k}}\right)^{1/k}
  \le
  \left(\frac{T_{j,G}(a)}{\binom{p}{j}}\right)^{1/j}.
\]
Equivalently,
\[
  T_{k,G}(a)
  \le
  \frac{\binom{p}{k}}{\binom{p}{j}^{k/j}}
  T_{j,G}(a)^{k/j}.
\]
Moreover, for each such pair $j,k$, the constant is sharp:
\[
  \sup_{\substack{a\in C^1(G)\\T_{j,G}(a)>0}}
  \frac{T_{k,G}(a)}{T_{j,G}(a)^{k/j}}
  =
  \frac{\binom{p}{k}}{\binom{p}{j}^{k/j}}.
\]

\end{maintheorem}

We begin with the Euclidean two-form inequality, which follows
from the skew-symmetric normal form and the classical Maclaurin
inequalities. We then prove a zero-propagation lemma: if all
order-$j$ clique Pfaffians containing a vertex vanish, then so
do all higher-order clique Pfaffians containing that vertex.
This allows us to delete vertices with zero local order-$j$
mass. For each fixed pair $j<k$, a deterministic adjustment
between nonadjacent active vertices then reduces the support
to a clique, keeping $T_{j,G}$ fixed without decreasing $T_{k,G}$.
Applying the Euclidean inequality to the resulting form proves
the graph inequality. Finally, equal edge values on a matching
inside a maximum clique establish sharpness.

\subsection{The Euclidean two-form inequality}

We first express the classical inequality in the exterior-algebra
notation of \cref{subsec:pfaffians-maclaurin}.

\begin{lemma}
\label{lem:two-form-maclaurin}
Let $\beta$ be a real two-form on an $N$-dimensional Euclidean
vector space, and put $r=\lfloor N/2\rfloor$. For all integers
$1\le j<k\le r$,
\[
  \left(
    \frac{\|\beta^k/k!\|^2}{\binom{r}{k}}
  \right)^{1/k}
  \le
  \left(
    \frac{\|\beta^j/j!\|^2}{\binom{r}{j}}
  \right)^{1/j}.
\]
\end{lemma}

\begin{proof}
By the orthogonal canonical form for real skew-symmetric
matrices \cite[Corollary~2.5.11(b)]{HornJohnson2012},
there is an orthonormal basis of covectors
$f^1,\ldots,f^N$ such that
\[
  \beta
  =\sum_{\alpha=1}^{r}
    \sigma_\alpha f^{2\alpha-1}\wedge f^{2\alpha},
\]
where some of the $\sigma_\alpha$ may be zero.
The two-form factors commute, and the square of each factor
$f^{2\alpha-1}\wedge f^{2\alpha}$ vanishes. Thus, for
$1\le s\le r$,
\[
  \frac{\beta^s}{s!}
  =\sum_{\alpha_1<\cdots<\alpha_s}
    \sigma_{\alpha_1}\cdots\sigma_{\alpha_s}
    f^{2\alpha_1-1}\wedge f^{2\alpha_1}
    \wedge\cdots\wedge
    f^{2\alpha_s-1}\wedge f^{2\alpha_s}.
\]
The exterior basis elements in this sum are orthonormal, so
\begin{equation}
  \left\|\frac{\beta^s}{s!}\right\|^2
  =e_s(\sigma_1^2,\ldots,\sigma_r^2).
  \label{eq:two-form-block-norm}
\end{equation}
Applying \cref{thm:classical-maclaurin} to the nonnegative
numbers $\sigma_1^2,\ldots,\sigma_r^2$ proves the inequality.
\end{proof}

For a complete graph, \eqref{eq:prelim-pfaffian-norm} identifies
$T_{s,G}(a)$ with $\|\beta_a^s/s!\|^2$, so the lemma applies
directly. For a general graph, $T_{s,G}$ includes only those
principal Pfaffians whose vertices form cliques. We next show
how to reduce this restricted sum to the complete-graph case
while retaining the clique-number normalization.

\subsection{Zero propagation}

For a vertex $v\in V(G)$ and an integer $j\ge1$, define
\[
  x_j(v;a)
  :=\sum_{\substack{S\in\Kcal_{2j}(G)\\v\in S}}
    \Pf\bigl(a[S]\bigr)^2.
\]
This is the contribution to $T_{j,G}(a)$ from cliques containing
$v$. The next lemma allows us to set all entries incident to
$v$ equal to zero whenever this contribution vanishes.

\begin{lemma}
\label{lem:zero-propagation}
Let $a\in C^1(G)$, $v\in V(G)$, and $j\ge1$. If
$x_j(v;a)=0$, then $x_k(v;a)=0$ for every integer $k\ge j$.
Consequently, setting all edge values incident to $v$ equal
to zero leaves $T_{k,G}(a)$ unchanged for every $k\ge j$.
\end{lemma}

\begin{proof}
Fix $k\ge j$ and a clique $T\in\Kcal_{2k}(G)$ containing $v$.
On $\mathbb R^T$, consider the two-form
\[
  \beta=\sum_{\substack{u<w\\u,w\in T}}
    a_{uw}\,e^u\wedge e^w.
\]
Every $2j$-subset of $T$ containing $v$ is a clique of $G$.
Since $x_j(v;a)=0$ is a sum of nonnegative terms, the
Pfaffian of each such subset is zero. By
\eqref{eq:prelim-pfaffian-exterior},
\[
  \frac{\beta^j}{j!}
  =\sum_{\substack{S\subseteq T\\|S|=2j}}
    \Pf\bigl(a[S]\bigr)e^S,
\]
so every coefficient involving $e^v$ vanishes.
Let $e_v$ be the coordinate vector dual to $e^v$, and let
$\iota_{e_v}$ denote contraction by this vector. We obtain
$\iota_{e_v}\beta^j=0$.

The contraction rule for a two-form gives
\[
  \iota_{e_v}(\beta^m)
  =m(\iota_{e_v}\beta)\wedge\beta^{m-1}
  \qquad(m\ge1).
\]
It follows that
\[
  \iota_{e_v}\beta^k
  =\frac{k}{j}
    \bigl(\iota_{e_v}\beta^j\bigr)\wedge\beta^{k-j}
  =0.
\]
On the other hand, $|T|=2k$ implies
\[
  \beta^k=k!\Pf\bigl(a[T]\bigr)e^T.
\]
Since $v\in T$, the form $\iota_{e_v}e^T$ is nonzero.
Thus $\Pf(a[T])=0$. This holds for every $2k$-clique
containing $v$, proving $x_k(v;a)=0$.

Setting the entries incident to $v$ equal to zero affects
only Pfaffians on vertex sets containing $v$. At every
order $k\ge j$, these Pfaffians were already zero and remain
zero, which proves the final assertion.
\end{proof}

\subsection{Reduction to the clique case}

A vertex $v$ is called \emph{active} for an edge form $a$ if
$a_{vw}\ne0$ for some $w$. We reduce the number of active
vertices by transferring their energy contributions between
nonadjacent vertices.

\begin{lemma}
\label{lem:transfer}
Fix integers $1\le j<k$. For every $a\in C^1(G)$, there exists
an edge form $b\in C^1(G)$ whose active vertices form a possibly
empty clique and such that
\[
  T_{j,G}(b)=T_{j,G}(a),
  \qquad
  T_{k,G}(b)\ge T_{k,G}(a).
\]
\end{lemma}

\begin{proof}
We successively modify the form, writing $a$ for the current
form at each step. First remove any active vertex $v$ with
$x_j(v;a)=0$ by setting all its incident entries equal to zero.
By \cref{lem:zero-propagation}, this preserves both $T_{j,G}$
and $T_{k,G}$. Repeat until every remaining active vertex has
positive local order-$j$ mass.

If the active vertices form a clique, stop. Otherwise, choose
two nonadjacent active vertices $u,v$ and put
\[
  x=x_j(u;a)>0,\qquad y=x_j(v;a)>0,
  \qquad X=x_k(u;a),\qquad Y=x_k(v;a).
\]
Interchanging $u$ and $v$ if necessary, assume $X/x\ge Y/y$.
Define $a'$ by setting all entries incident to $v$ equal to
zero and multiplying all entries incident to $u$ by
$\sqrt{(x+y)/x}$. All other entries remain unchanged.
Since $u,v$ are nonadjacent, these operations do not interfere,
and $a'$ is still an edge form on $G$.

No clique contains both $u$ and $v$, so we may write
\[
  T_{j,G}(a)=C_j+x+y,
  \qquad T_{k,G}(a)=C_k+X+Y,
\]
where $C_j,C_k$ are the contributions from cliques containing
neither vertex.
Every perfect matching of a clique containing $u$ uses
exactly one edge incident to $u$. Thus the adjustment multiplies
its squared Pfaffian by $(x+y)/x$, while all Pfaffians on
cliques containing $v$ become zero. Consequently,
\[
  T_{j,G}(a')=C_j+\frac{x+y}{x}x=T_{j,G}(a),
\]
and
\[
  T_{k,G}(a')-T_{k,G}(a)
  =\frac{y}{x}X-Y
  =y\left(\frac{X}{x}-\frac{Y}{y}\right)\ge0.
\]

Repeat, removing vertices with zero local order-$j$ mass
before each adjustment. Every removal or adjustment strictly
decreases the number of active vertices and creates no new
ones. The process therefore stops after at most $|V(G)|$
operations. Its final form $b$ has clique support and satisfies
the two asserted energy relations. This also covers the zero
form, whose active set is empty.
\end{proof}

\begin{remark}
This is an adjustment argument: each step preserves the
lower-order energy and does not decrease the chosen higher-order
energy, while reducing the number of active vertices. The
resulting form $b$ may depend on both $j$ and $k$.
\end{remark}

\subsection{Proof and sharpness}

\begin{proof}[Proof of \cref{thm:masked-maclaurin}]
Suppose first that $T_{j,G}(a)=0$ for some integer $j\ge1$.
Then $x_j(v;a)=0$ for every vertex $v$.
By \cref{lem:zero-propagation}, $x_k(v;a)=0$ for every
$v$ and every $k\ge j$. Every $2k$-clique therefore has zero
Pfaffian, and $T_{k,G}(a)=0$.

Now fix $1\le j<k\le p$ and assume $T_{j,G}(a)>0$.
Apply \cref{lem:transfer} to obtain an edge form $b$ with
\[
  T_{j,G}(b)=T_{j,G}(a),
  \qquad
  T_{k,G}(b)\ge T_{k,G}(a),
\]
whose active vertices form a clique $C$. Associate to $b$
the two-form
\[
  \beta_b=\sum_{u<v}b_{uv}\,e^u\wedge e^v.
\]
Any principal submatrix involving a vertex outside $C$ has a
zero row and column, so its Pfaffian vanishes. Every subset
of $C$ is a clique. Consequently,
\eqref{eq:prelim-pfaffian-norm} gives
\[
  T_{s,G}(b)=\left\|\frac{\beta_b^s}{s!}\right\|^2
  \qquad(s\ge1).
\]

Since $|C|\le\omega(G)$, the skew normal form of $\beta_b$
has at most $\lfloor|C|/2\rfloor\le p$ nonzero blocks.
Append zeros to their parameters to obtain a list
$\sigma_1,\ldots,\sigma_p$. By the calculation in
\eqref{eq:two-form-block-norm},
\[
  T_{s,G}(b)=e_s(\sigma_1^2,\ldots,\sigma_p^2).
\]
Applying \cref{thm:classical-maclaurin} to this list gives
\[
  T_{k,G}(a)\le T_{k,G}(b)
  \le\frac{\binom{p}{k}}{\binom{p}{j}^{k/j}}
     T_{j,G}(b)^{k/j}
  =\frac{\binom{p}{k}}{\binom{p}{j}^{k/j}}
     T_{j,G}(a)^{k/j}.
\]

To prove sharpness, choose a $2p$-clique of $G$ and a perfect
matching of its vertices. Give each matching edge the same
nonzero value $\tau$ in the increasing orientation, and set
all other entries equal to zero. A principal Pfaffian of
order $2s$ is nonzero precisely when its vertex set is the
union of $s$ of these matching edges. There are
$\binom{p}{s}$ such sets, all of which are cliques, and each
has squared Pfaffian $\tau^{2s}$. Hence
\[
  T_{s,G}(a)=\binom{p}{s}\tau^{2s},
  \qquad
  \left(\frac{T_{s,G}(a)}{\binom{p}{s}}\right)^{1/s}
  =\tau^2
  \qquad(1\le s\le p).
\]
Equality therefore holds simultaneously for every admissible
pair $j,k$, proving the sharpness assertion.
\end{proof}

The case $j=1$ and $k=2$ bounds $T_{2,G}(a)$ in terms of
$E_G(a)^2$.

\begin{corollary}
\label{cor:graph-four-point}
Let $p=\lfloor\omega(G)/2\rfloor\ge1$. Then 
for every $a\in C^1(G)$,
\[
  T_{2,G}(a)\le\frac{p-1}{2p}E_G(a)^2.
\]

\end{corollary}

\begin{proof}
If $p=1$, then $\omega(G)\le3$, so $G$ has no four-vertex
clique and $T_{2,G}(a)=0$. If $p\ge2$, apply
\cref{thm:masked-maclaurin} with $j=1$ and $k=2$, using
$T_{1,G}(a)=E_G(a)$:
\[
  T_{2,G}(a)
  \le\frac{\binom{p}{2}}{\binom{p}{1}^2}E_G(a)^2
  =\frac{p-1}{2p}E_G(a)^2.
\]

\end{proof}

Define
\begin{equation}
  D_G(a):=E_G(a)^2-2T_{2,G}(a).
  \label{eq:clique-certificate}
\end{equation}
If $a\ne0$, then $G$ has an edge, so $p\ge1$ and $E_G(a)>0$.
By \cref{cor:graph-four-point},
\[
  D_G(a)\ge\frac{E_G(a)^2}{p}>0.
\]
Consequently,
\begin{equation}
  0<D_G(a),\qquad
  \frac{2E_G(a)^2}{D_G(a)}
  \le2\left\lfloor\frac{\omega(G)}2\right\rfloor
  \le\omega(G).
  \label{eq:certificate-introduction}
\end{equation}
Thus a large value of this quotient certifies a large clique.
In \cref{sg:section}, we construct edge forms with a large
quotient from spectral information and use
\eqref{eq:certificate-introduction} to prove
\cref{sg:theorem}.

\section{The spectral clique bound}
\label{sg:section}
This section concerns how large a clique must occur in a regular graph
whose smallest adjacency eigenvalue is at least $-3$. For graphs
satisfying this spectral condition and having sufficiently large minimum
degree, Koolen, Yang and Yang~\cite{KYY} developed a structural
description in terms of Hoffman graphs. Yang and
Koolen~\cite[Corollary~1.5]{YangKoolen2025} subsequently showed that,
for sufficiently large $d$, every vertex of a connected $d$-regular
graph in this class lies in a clique of order at least
$(d-\varepsilon)/15$, where $\varepsilon>0$ is an absolute constant.
In particular, this gives $\omega(G)\ge d/15-O(1)$. We combine
the KYY structure with the helicity energy certificate
\eqref{eq:certificate-introduction} to improve this global bound to
$\omega(G)\ge d/4-O(1)$, with a rank-dependent refinement.

The starred lemmas below have elementary but lengthy combinatorial
counting proofs, which may be skipped on a first reading to focus
on the overall proof strategy.
We prove the following rank-dependent bound:

\begin{maintheorem}\label{sg:theorem}
There is an absolute integer $B_*\ge1$ such that the following holds.
Let $G$ be a nonempty simple $d$-regular graph on $n$ vertices, let $M_G$ be
its adjacency matrix, and assume $\lambda_{\min}(M_G)\ge-3$. Put
\[
 \eta=\frac{\rank(M_G+3I)}{n},\qquad C_*=2B_*+6,\qquad
 g(\eta)=
 \begin{cases}
 \dfrac{1-2\eta}{3-4\eta},&0\le\eta\le\tfrac14,\\[4pt]
 \dfrac{2\eta+1}{6},&\tfrac14\le\eta\le\tfrac12,\\[4pt]
 \dfrac13,&\tfrac12\le\eta\le1.
 \end{cases}
\]
Then
\[
                         \omega(G)\ge d\,g(\eta)-C_*.
\]
More precisely, if $d\ge6B_*$, there is a nonzero edge form
$a\in C^1(G)$ for which
\begin{equation}\label{sg:certificate-target}
 \frac{2E_G(a)^2}{D_G(a)}\ge d\,g(\eta)-C_*,
 \qquad D_G(a):=E_G(a)^2-2T_{2,G}(a)>0.
\end{equation}
The form is constructed on an induced subgraph of $G$ and extended
by zero to all of $G$.
\end{maintheorem}

The clique bound in Theorem~\ref{sg:theorem} is immediate when
$d<6B_*$, since $g(\eta)\le1/3$ and $C_*=2B_*+6$ give
$d\,g(\eta)-C_*<0$. Throughout the remainder of this section,
we therefore assume $d\ge6B_*$, with $B_*$ chosen below.

The proof constructs a form on a suitable fat
neighborhood in the associated Hoffman graph defined below.
If some vertex has exactly one fat neighbor, a local spectral-radius
bound supplies the required form. Otherwise, residual classes and
weighted incidence counts bound the number of components of a local
complement, and we use its Laplacian kernel. In both cases, zero
extension and Corollary~\ref{cor:graph-four-point} turn the local
energy bound into the global clique bound.

The proof is organized into four parts. We first introduce the
associated Hoffman graph and establish the coverage and local
spectral bounds. We then give two constructions of edge forms,
using either the connected components or the spectral radius of
a complement graph. Next, we analyze the residual Gram vectors
and derive the weighted incidence estimates needed to relate
local structure to the rank of $M_G+3I$. Finally, we obtain the
rank bounds, select a suitable local form in each case, and
complete the proof of Theorem~\ref{sg:theorem}.

\subsection{Hoffman graphs, coverage, and local spectral bounds}\label{subsec:hoffman}

We retain the graph and notation of Theorem~\ref{sg:theorem}.

A \emph{Hoffman graph} is a graph
$H=(V_s,V_f,E)$ whose vertices are partitioned into \emph{slim} and
\emph{fat} vertices, such that the fat vertices are pairwise nonadjacent
and each fat vertex has a slim neighbor. Its \emph{slim graph} is the
graph induced on $V_s$. When that slim graph is $G$, define the slim--fat incidence matrix
$C=(C_{vf})_{v\in V_s,\ f\in V_f}$ by
\[
C_{vf}=
\begin{cases}
1,&\text{if }vf\in E(H),\\
0,&\text{otherwise}.
\end{cases}
\]
Then the adjacency matrix of $H$ is

\[
\begin{pmatrix}M_G&C\\ C^T&0\end{pmatrix}.
\]

The \emph{special matrix} of $H$ is
\[
 B(H):=M_G-CC^T.
\]
Following \cite[Definition~2.4]{KYY}, we define the
\emph{eigenvalues of the Hoffman graph} $H$ to be the
eigenvalues of its special matrix $B(H)$. In particular,
\[
  \lambda_{\min}(H):=\lambda_{\min}(B(H)).
\]
This refers to the spectrum of the special matrix, rather
than that of the adjacency matrix of the underlying graph.
When $H$ has no fat vertices, $B(H)=M_G$, so this definition
recovers the usual graph spectrum.

The associated-graph construction requires the following exclusion,
which follows from \cref{thm:prelim-interlacing} and a Rayleigh
quotient calculation.
\begin{lemma}
Let $\widetilde K_{24}$ be the graph obtained from a clique on
$24$ vertices by adjoining one vertex adjacent to exactly $12$
clique vertices. If a finite simple graph $G$ satisfies
$\lambda_{\min}(M_G)\ge-3$, then $G$ contains no induced subgraph
isomorphic to $\widetilde K_{24}$.
\end{lemma}

\begin{proof}
Since $\widetilde K_{24}$ has $25$ vertices, we construct a vector
\[
  x=(x_u)_{u\in V(\widetilde K_{24})}
  \in\R^{V(\widetilde K_{24})}\cong\R^{25},
\]
whose coordinates are indexed by the vertices in the same order
as the rows and columns of $M_{\widetilde K_{24}}$.
Let $v$ denote the additional vertex, let $A$ be its set of twelve
neighbors in the clique, and let $B$ be the set of the remaining
twelve clique vertices. Thus
$V(\widetilde K_{24})=\{v\}\sqcup A\sqcup B$.
Choose one vertex $b_0\in B$ and define
\[
  x_u=
  \begin{cases}
    -4,&u=v,\\
     1,&u\in A,\\
    -1,&u\in B\setminus\{b_0\},\\
     0,&u=b_0.
  \end{cases}
\]
The resulting vector satisfies
\[
  \|x\|^2=\sum_{u\in V(\widetilde K_{24})}x_u^2
  =16+12+11=39.
\]
Put $K=A\cup B$. Each unordered edge $\{u,w\}$ contributes
$2x_ux_w$ to the adjacency quadratic form, since both matrix
entries $(u,w)$ and $(w,u)$ are included. All pairs of distinct
vertices in $K$ are edges, and the edges incident to $v$ are
precisely those joining $v$ to $A$. Hence
\[
  x^TM_{\widetilde K_{24}}x
  =2\sum_{\substack{u<w\\u,w\in K}}x_ux_w
   +2x_v\sum_{u\in A}x_u.
\]
Expanding the square of the sum of the clique coordinates gives
\[
  2\sum_{\substack{u<w\\u,w\in K}}x_ux_w
  =\left(\sum_{u\in K}x_u\right)^2-\sum_{u\in K}x_u^2.
\]
By the definition of $x$,
\[
  \sum_{u\in K}x_u=12-11+0=1,\qquad
  \sum_{u\in K}x_u^2=12+11=23,\qquad
  \sum_{u\in A}x_u=12.
\]
Substituting these values and $x_v=-4$, we obtain
\[
  \begin{aligned}
    x^TM_{\widetilde K_{24}}x
    &=1^2-23+2(-4)(12)\\
    &=-22-96=-118<-117=-3\|x\|^2.
  \end{aligned}
\]
The Rayleigh quotient therefore gives
\[
\lambda_{\min}(M_{\widetilde K_{24}})
\le-\frac{118}{39}<-3.
\]
If $G$ contained an induced copy of $\widetilde K_{24}$, eigenvalue
interlacing would imply
\[
\lambda_{\min}(M_G)
\le\lambda_{\min}(M_{\widetilde K_{24}})<-3,
\]
contrary to the hypothesis.
\end{proof}

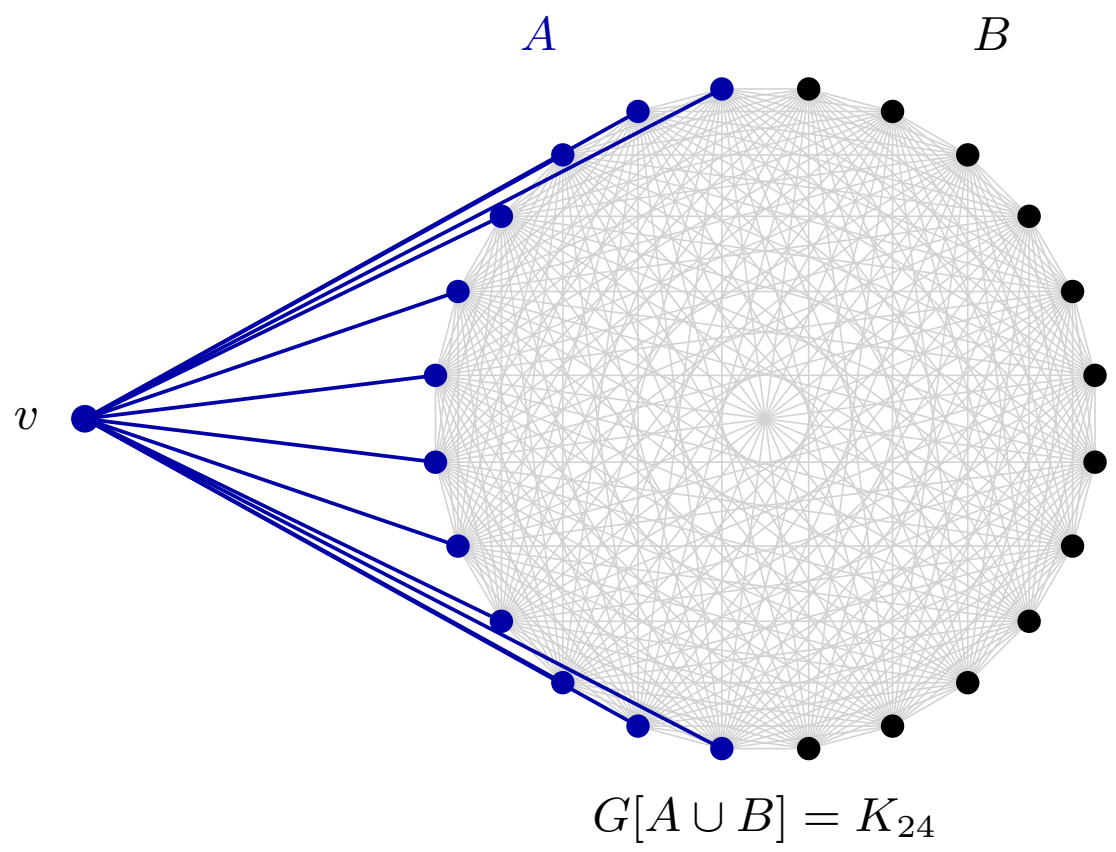
\begin{figure}[htbp]
  \centering
  \begin{tikzpicture}[
    line cap=round,
    line join=round,
    font=\small
  ]
    \foreach \i in {1,...,24} {
      \coordinate (c\i) at ({97.5+15*(\i-1)}:2.2);
    }

    \foreach \i in {1,...,23} {
      \pgfmathtruncatemacro{\nextvertex}{\i+1}
      \foreach \j in {\nextvertex,...,24} {
        \draw[gray!35,line width=0.3pt]
          (c\i) -- (c\j);
      }
    }

    \coordinate (v) at (-4.5,0);
    \foreach \i in {1,...,12} {
      \draw[blue!65!black,line width=0.7pt]
        (v) -- (c\i);
    }

    \foreach \i in {1,...,12} {
      \fill[blue!65!black] (c\i) circle (2.2pt);
    }

    \foreach \i in {13,...,24} {
      \fill[black] (c\i) circle (2.2pt);
    }

    \fill[blue!65!black] (v) circle (2.6pt);
    \node[left=5pt] at (v) {$v$};

    \node[blue!65!black] at (-1.5,2.55) {$A$};
    \node at (1.5,2.55) {$B$};
    \node at (0,-2.65) {$G[A\cup B]=K_{24}$};
  \end{tikzpicture}
  \caption{The graph $\widetilde K_{24}$.
  The sets $A$ and $B$ each contain twelve vertices, and their
  union induces a complete graph, whose edges are shown in gray.
  The additional vertex $v$ is adjacent to every vertex of $A$
  and to no vertex of $B$; its twelve incident edges are shown
  in blue.}
  \label{fig:forbidden-k24}
\end{figure}

We use the associated-graph construction of \cite[Section~2.5]{KYY}
with parameter $m=12$. 

For an integer $r\ge169$, let $\mathcal C_r$ be the family
of maximal cliques of $G$ of order at least $r$. For any clique $K$ put
\[
 b_K(v)=\bigl|\{u\in K\setminus\{v\}:uv\notin E(G)\}\bigr|,
 \qquad Q(K)=\{v\in V(G):b_K(v)\le11\}.
\]
On $\mathcal C_r$ define
\[
 K\sim L
 \quad\Longleftrightarrow\quad
 b_L(v)\le11\ \text{for every }v\in K,
 \quad
 b_K(u)\le11\ \text{for every }u\in L.
\]
Section~2.5 of \cite{KYY} states that, in the absence of an induced
$\widetilde K_{24}$, this is an equivalence relation and $Q(K)$
depends only on the class $[K]$. The \emph{associated Hoffman graph}
$\mathfrak g(G,12,r)$ retains $G$ as its slim graph and has exactly
one fat vertex $f_{[K]}$ for each equivalence class, with slim
neighborhood $Q(K)$. In particular,
\begin{equation}\label{sg:clique-covered}
 K\subseteq Q(K)=Q(L)=N(f_{[K]})
 \qquad K\sim L.
\end{equation}

Thus every large maximal clique is contained in the fat neighborhood
assigned to its equivalence class, independently of the chosen representative.

We apply the following result of \cite[Lemma~5.1]{KYY}.
\begin{lemma}[Koolen--Yang--Yang]\label{sg:kyy}
There is an absolute positive integer $r_0$ such that, for every
connected graph $G$ with $\lambda_{\min}(M_G)\ge-3$ and every integer
$r\ge\max\{r_0,169\}$, the associated Hoffman graph
$\mathfrak g(G,12,r)$ has smallest eigenvalue at least $-3$.
Consequently, its slim--fat incidence matrix satisfies
\begin{equation}\label{sg:residual-psd}
 R:=M_G-CC^T+3I\succeq0.
\end{equation}
\end{lemma}

For a disconnected graph, we apply the associated-graph construction to each connected component and take the disjoint union of the resulting Hoffman graphs. The slim--fat incidence matrix and the residual matrix are then block diagonal across components, so \eqref{sg:residual-psd} holds in this case as well. All counts below are taken over the entire graph.

Fix once and for all an integer $r_*\ge\max\{r_0,169\}$, and set
\[
 B_*:=\mathrm{R}(r_*-1,10)-1,
 \qquad H:=\mathfrak g(G,12,r_*).
\]
Here $\mathrm{R}(s,t)$ is the Ramsey number defined in
\cref{subsec:spectral-preliminaries}.

\begin{definition}
For a slim vertex $v\in V(G)$, let
\[
F(v)=\{f\in V_f:vf\in E(H)\}
\]
be its set of fat neighbors. For a fat vertex $f\in V_f$, define
its \emph{slim neighborhood} and its order by
\[
Q_f=\{v\in V(G):f\in F(v)\},
\qquad q_f=|Q_f|,
\]
and write
\[
H_f=\overline{G[Q_f]}
\]
for the complement of the graph induced by this neighborhood.

An edge $uv\in E(G)$ is called \emph{uncovered} if
$F(u)\cap F(v)=\varnothing$. For each $v\in V(G)$, let
\[
N_u(v)=\{u\in N_G(v):F(u)\cap F(v)=\varnothing\}
\]
be the set of its uncovered neighbors.
\end{definition}

\begin{figure}[htbp]
  \centering
  \begin{tikzpicture}[
    line width=0.7pt,
    line cap=round,
    line join=round,
    font=\small,
    slim/.style={
      circle,draw=blue!65!black,fill=blue!8,
      minimum size=4.8mm,inner sep=0pt
    },
    fat/.style={
      rectangle,draw=blue!65!black,fill=blue!18,
      minimum size=5.5mm,inner sep=0pt
    },
    incidence/.style={blue!65!black,dashed}
  ]
    \begin{scope}
      \fill[blue!4,rounded corners=5pt]
        (-0.92,-1.3) rectangle (1.32,0.72);
      \coordinate (h1) at (-0.6,0.4);
      \coordinate (h2) at (1,0.4);
      \coordinate (h3) at (-0.6,-1);
      \coordinate (h4) at (1,-1);
      \coordinate (h5) at (-1.8,0.4);
      \coordinate (hf) at (0.55,2);
      \coordinate (hg) at (-1.3,2);

      \foreach \i/\j in {1/3,1/4,2/3,2/4,3/4,1/5}
        \draw (h\i)--(h\j);
      \foreach \i in {1,...,4}
        \draw[incidence] (hf)--(h\i);
      \foreach \i in {1,5}
        \draw[incidence] (hg)--(h\i);

      \foreach \i in {1,...,4}
        \node[slim] at (h\i) {$\i$};
      \node[slim,draw=gray,fill=gray!8] at (h5) {$5$};
      \node[fat] at (hf) {$f$};
      \node[fat] at (hg) {$g$};
      \node[blue!65!black] at (0.2,-1.55) {$Q_f$};
      \node at (-0.25,-2.1) {(a) Hoffman graph};
    \end{scope}

    \node at (2.5,-0.25) {$\longmapsto$};

    \begin{scope}[xshift=4.6cm]
      \coordinate (q1) at (-0.8,0.4);
      \coordinate (q2) at (0.8,0.4);
      \coordinate (q3) at (-0.8,-1);
      \coordinate (q4) at (0.8,-1);
      \foreach \i/\j in {1/3,1/4,2/3,2/4,3/4}
        \draw (q\i)--(q\j);
      \foreach \i in {1,...,4}
        \node[slim] at (q\i) {$\i$};
      \node at (0,-2.1) {(b) $G[Q_f]$};
    \end{scope}

    \node at (6.85,-0.25) {$\longmapsto$};

    \begin{scope}[xshift=9.1cm]
      \coordinate (c1) at (-0.8,0.4);
      \coordinate (c2) at (0.8,0.4);
      \coordinate (c3) at (-0.8,-1);
      \coordinate (c4) at (0.8,-1);
      \draw[blue!65!black] (c1)--(c2);
      \foreach \i in {1,...,4}
        \node[slim] at (c\i) {$\i$};
      \node at (0,-2.1) {(c) $H_f=\overline{G[Q_f]}$};
    \end{scope}
  \end{tikzpicture}
  \caption{The notation illustrated on a small Hoffman graph.
  Squares are fat vertices, circles are slim vertices, dashed
  lines are slim--fat edges, and solid lines join slim vertices.
  The shaded set $Q_f$ consists of the four slim neighbors of $f$.
  Taking the induced slim graph and then its complement leaves
  only the edge $12$; the vertices $3,4$ remain as isolated vertices.}
  \label{fig:fat-neighborhood-complement}
\end{figure}
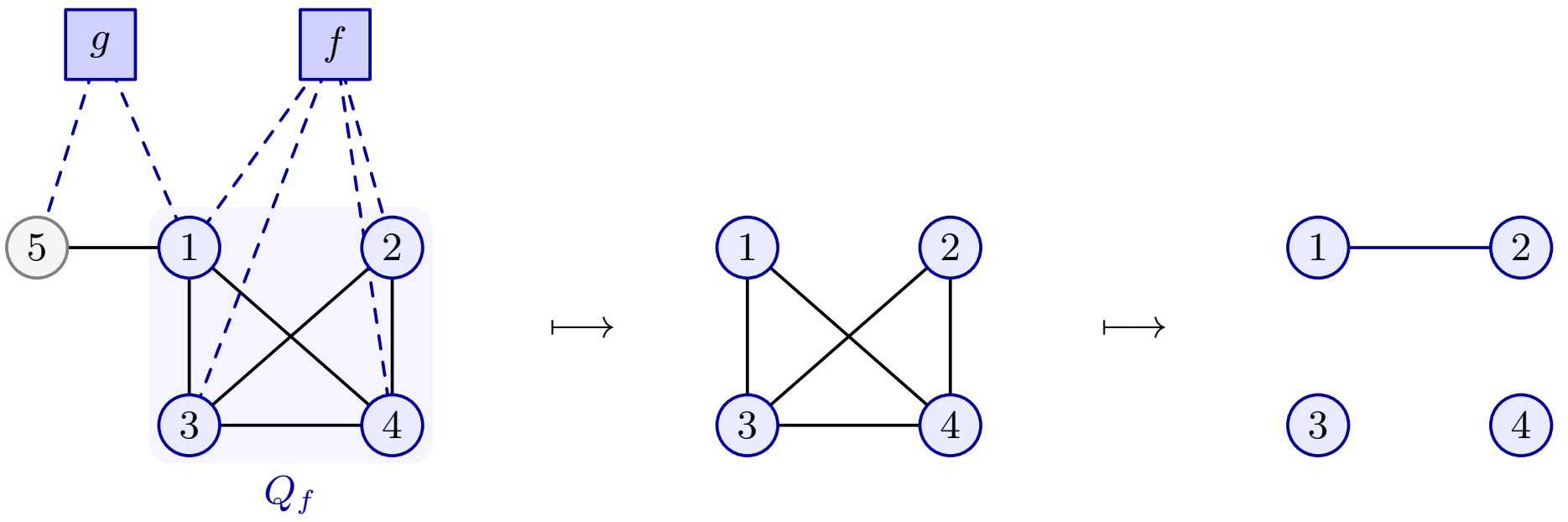

To illustrate these definitions independently of the size assumptions
above, consider the small Hoffman graph in
\cref{fig:fat-neighborhood-complement}. Here
\[
  F(1)=\{f,g\},\qquad
  Q_f=\{1,2,3,4\},\qquad q_f=4.
\]
The graph $G[Q_f]$ contains every edge on $Q_f$ except $12$.
Consequently, $H_f$ has the same four slim vertices, with $12$
as its only edge and $3,4$ isolated. The fat vertex $f$ serves
only to specify the set $Q_f$.

\begin{lemma}\label{sg:coverage}
Every vertex $v$ satisfies $|N_u(v)|\le B_*$ and
$1\le |F(v)|\le3$. With
$t=\min_{v\in V(G)}|F(v)|$, every fat neighborhood satisfies
\[
                         \rho(H_f)\le3-t.
\]
\end{lemma}
\begin{proof}
By \eqref{eq:prelim-no-large-star}, $G$ has no induced
$K_{1,10}$. Suppose that
$|N_u(v)|\ge B_*+1=\mathrm{R}(r_*-1,10)$.
An independent ten-set in $G[N_u(v)]$, together with $v$, would
induce $K_{1,10}$. By \cref{thm:prelim-ramsey}, there is a clique
$S\subseteq N_u(v)$ of order $r_*-1$. Since each vertex of $S$
is adjacent to $v$, the set $S\cup\{v\}$ is an $r_*$-clique.
Extend it to a maximal clique $K\in\mathcal C_{r_*}$.
By \eqref{sg:clique-covered}, all vertices of $S\cup\{v\}$ belong
to $Q_{f_{[K]}}$. Thus
$f_{[K]}\in F(u)\cap F(v)$ for every $u\in S$, contrary to
$S\subseteq N_u(v)$. This proves $|N_u(v)|\le B_*$.

If $F(v)=\varnothing$, then all $d$ neighbors of $v$
are uncovered, a contradiction. Thus $|F(v)|\ge1$. On the other
hand, the nonnegative diagonal of \eqref{sg:residual-psd} gives
\[
                  0\le R_{vv}=3-|F(v)|,
\]
so $|F(v)|\le3$ and $t\in\{1,2,3\}$.

Fix $f\in V_f$. For distinct $u,v\in Q_f$ we have
$|F(u)\cap F(v)|\ge1$, and hence
\[
 R_{uv}=(M_G)_{uv}-|F(u)\cap F(v)|
 \le (M_G)_{uv}-1=-(M_{H_f})_{uv},
 \qquad R_{vv}\le3-t.
\]
By \cref{thm:prelim-perron}, choose a nonzero nonnegative
eigenvector $x$ of $M_{H_f}$ with eigenvalue $\rho(H_f)$. The principal submatrix $R[Q_f]$ is positive
semidefinite, so entrywise comparison against $x\ge0$ yields
\[
 0\le x^TR[Q_f]x
 \le (3-t)\|x\|^2-x^TM_{H_f}x
 =\bigl(3-t-\rho(H_f)\bigr)\|x\|^2.
\]
Since $x\ne0$, this proves $\rho(H_f)\le3-t$.
\end{proof}

\begin{lemma}
Let $t=\min_{v\in V(G)}|F(v)|$.
For every vertex $v$ with $|F(v)|=t$, there exists a fat neighbor
$f\in F(v)$ such that
\[
q_f\ge 1+\frac{d-B_*}{t}.
\]
\end{lemma}

\begin{proof}
By Lemma~\ref{sg:coverage}, we have $|N_u(v)|\le B_*$ and $t\ge1$.
Every covered neighbor of $v$ belongs to $Q_f\setminus\{v\}$
for some $f\in F(v)$. Therefore
\begin{equation}\label{sg:fat-average}
d-B_*
\le |N_G(v)\setminus N_u(v)|
\le \sum_{f\in F(v)}(q_f-1).
\end{equation}
There are exactly $t$ summands, so at least one is at least
$(d-B_*)/t$. This proves the assertion.
\end{proof}

\subsection{Two local edge-form constructions}

We construct edge forms from the components or the spectral radius
of a complement graph. For an edge form on a graph $X$, recall
$D_X(a)=E_X(a)^2-2T_{2,X}(a)$. We first isolate a matrix construction
used in both arguments. It reduces the support and energy
calculations to conditions on an orthonormal family of vectors.

\begin{lemma}
\label{sg:orthonormal-construction}
Let $X$ be a finite simple graph on a vertex set $Q$, and let
$p\ge0$ be an integer with $2p\le |Q|$. Suppose that
\[
  U=(u_1\ \cdots\ u_{2p})\in\R^{|Q|\times 2p},
  \qquad U^TU=I_{2p}.
\]
For each $x\in Q$, let $r_x$ denote the row of $U$ indexed by $x$.
Assume that $r_x$ and $r_y$ are linearly dependent whenever
$x\ne y$ and $\{x,y\}\notin E(X)$. Set
\[
  \Omega_p=\operatorname{diag}\left(
    \underbrace{
      \begin{pmatrix}0&1\\-1&0\end{pmatrix},\ldots,
      \begin{pmatrix}0&1\\-1&0\end{pmatrix}
    }_{p\text{ blocks}}
  \right),
  \qquad a=U\Omega_pU^T.
\]
Then $a$ is an edge form on $X$, and
$\Pf(a[S])=0$ for every even vertex set $S$ that is not a clique
of $X$. Moreover,
\begin{equation}
  T_{j,X}(a)=\binom pj\quad(j\ge1),
  \qquad E_X(a)=D_X(a)=p,
  \label{sg:orthonormal-energies}
\end{equation}
where $\binom pj=0$ for $j>p$. If $p>0$, then
\[
  \frac{2E_X(a)^2}{D_X(a)}=2p.
\]
For $p=0$, the construction is interpreted as the zero form.
\end{lemma}

\begin{proof}
The case $p=0$ is immediate, so assume $p>0$.
The matrix $a$ is skew-symmetric. If $x$ and $y$ are distinct
nonadjacent vertices of $X$, linear dependence of $r_x,r_y$ and
skew-symmetry of $\Omega_p$ give
\[
  a_{xy}=r_x\Omega_pr_y^T=0.
\]
Thus $a$ vanishes on every nonedge of $X$.

Now let $S$ be an even vertex set that is not a clique, and choose
nonadjacent $x,y\in S$. The corresponding rows of $a[S]$ are
linearly dependent, since they are obtained by multiplying
$r_x$ and $r_y$ by the same matrix. Hence $a[S]$ is singular and
\[
  \Pf(a[S])^2=\det(a[S])=0.
\]

To compute the energies, identify each column $u_i$ with the
covector
\[
  u_i^\flat=\sum_{x\in Q}U_{xi}e^x,
\]
where $(e^x)_{x\in Q}$ is the standard orthonormal basis of
covectors. Each $2\times2$ block of $\Omega_p$ contributes
$u_{2r-1}u_{2r}^T-u_{2r}u_{2r-1}^T$ to the matrix product, so
\[
  a=U\Omega_pU^T
   =\sum_{r=1}^{p}
      \bigl(u_{2r-1}u_{2r}^T-u_{2r}u_{2r-1}^T\bigr).
\]
In particular, its $(x,y)$ entry is
\[
  a_{xy}=\sum_{r=1}^{p}
    \bigl(U_{x,2r-1}U_{y,2r}-U_{x,2r}U_{y,2r-1}\bigr).
\]
On the other hand, expanding the wedge product gives
\[
  \begin{aligned}
    u_{2r-1}^\flat\wedge u_{2r}^\flat
    &=\sum_{x,y\in Q}U_{x,2r-1}U_{y,2r}\,e^x\wedge e^y\\
    &=\sum_{x<y}
      \bigl(U_{x,2r-1}U_{y,2r}-U_{x,2r}U_{y,2r-1}\bigr)
      e^x\wedge e^y.
  \end{aligned}
\]
Here the diagonal terms vanish because $e^x\wedge e^x=0$,
and the terms indexed by $(x,y)$ and $(y,x)$ combine using
$e^y\wedge e^x=-e^x\wedge e^y$. Summing over $r$, the
coefficient of $e^x\wedge e^y$ is therefore exactly $a_{xy}$.
Thus the two-form associated with $a$ is
\[
  \beta=\sum_{x<y}a_{xy}e^x\wedge e^y
       =\sum_{r=1}^{p}u_{2r-1}^\flat\wedge u_{2r}^\flat.
\]
For $1\le j\le p$, its normalized exterior power is
\[
  \frac{\beta^j}{j!}
  =\sum_{1\le r_1<\cdots<r_j\le p}
    u_{2r_1-1}^\flat\wedge u_{2r_1}^\flat\wedge\cdots
    \wedge u_{2r_j-1}^\flat\wedge u_{2r_j}^\flat.
\]
The covectors $u_1^\flat,\ldots,u_{2p}^\flat$ are orthonormal,
so the $\binom pj$ summands are mutually orthogonal unit forms.
For $j>p$, the exterior power vanishes. Since every nonclique
Pfaffian vanishes, \eqref{eq:prelim-pfaffian-norm} now gives
\[
  T_{j,X}(a)
  =\sum_{\substack{S\subseteq Q\\|S|=2j}}\Pf(a[S])^2
  =\left\|\frac{\beta^j}{j!}\right\|^2
  =\binom pj.
\]
In particular, $E_X(a)=p$ and
\[
  D_X(a)=p^2-2\binom p2=p.
\]
The claimed quotient follows when $p>0$.
\end{proof}

To apply the lemma, it suffices to choose orthonormal columns so
that the rows at each pair of nonadjacent vertices are linearly
dependent. In the applications below, these rows are either equal
or one of them is zero. We first use normalized indicators of
components of the complement.

\begin{lemma}\label{sg:kernel}
Let $X$ be a finite simple graph on a vertex set $Q$, let
$J=\overline X$, and put
\[
 p_0=\left\lfloor\frac{c(J)}2\right\rfloor,
\]
where $c(J)$ is the number of connected components of $J$.
There is an explicit edge form $a$ on $X$ such that
\[
 E_X(a)=D_X(a)=p_0,
 \qquad T_{j,X}(a)=\binom{p_0}{j}\quad(j\ge1),
\]
with $\binom{p_0}{j}=0$ when $j>p_0$. If $p_0>0$, then
\[
 \frac{2E_X(a)^2}{D_X(a)}=2p_0\ge c(J)-1.
\]
If $p_0=0$, the construction gives the zero form and no quotient is taken.
\end{lemma}

\begin{proof}
If $p_0=0$, take $a=0$. Otherwise, choose $2p_0$ distinct
components $W_1,\ldots,W_{2p_0}$ of $J$, and put
\[
  u_i=\frac{\mathbf 1_{W_i}}{\sqrt{|W_i|}},
  \qquad U_0=(u_1\ \cdots\ u_{2p_0}).
\]
These columns are orthonormal. They also belong to $\ker L_J$,
since the Laplacian identity \eqref{eq:prelim-laplacian} shows
that this kernel consists of vectors constant on each component.

If $x,y$ are nonadjacent in $X$, then $\{x,y\}\in E(J)$,
so they lie in the same component of $J$ and have equal rows
in $U_0$. Lemma~\ref{sg:orthonormal-construction} therefore
applies to
\[
  a=U_0\Omega_{p_0}U_0^T
\]
and gives all the stated energy identities. When $p_0>0$,
\[
  \frac{2E_X(a)^2}{D_X(a)}
  =2p_0
  =2\left\lfloor\frac{c(J)}2\right\rfloor
  \ge c(J)-1.
\]
\end{proof}

\begin{remark}
If only a bound on the clique number is required, choosing one
vertex from each connected component of $J$ directly gives a
clique of order $c(J)$ in $X=\overline J$.
Lemma~\ref{sg:kernel} demonstrates how the helicity energy
framework translates this component structure into an explicit
edge form with exactly computable energies and, when $c(J)\ge2$,
a clique certificate.
\end{remark}

Although the next lemma introduces $J$ before $X$, in both
lemmas $J=\overline X$ and the edge form is constructed on $X$.

\begin{lemma}
\label{sg:radius}
Let $J$ be a finite simple graph on $N\ge4$ vertices with spectral
radius $\rho(J)\le2$, and let $X=\overline J$. There is a nonzero
edge form $a$ on $X$ such that
\[
 E_X(a)=D_X(a)=S>0,
 \qquad 2S\ge \frac N3-1.
\]
\end{lemma}

\begin{proof}
We will construct forms on the components of $J$, except that four types
of small components will be treated together. This common exceptional
block will incur a loss of at most one in the final bound.

By Smith's classification (\cref{thm:smith-classification}),
every connected component of $J$ is either a cycle or a tree
of one of the types
\[
A_n,\quad D_n,\quad E_6,\ E_7,\ E_8,\quad
\widetilde D_n,\quad
\widetilde E_6,\ \widetilde E_7,\ \widetilde E_8.
\]
The descriptions and diagrams in
\cref{subsec:spectral-preliminaries} show that each tree contains a path
omitting at most two vertices. Since every path in a tree is induced,
a tree component of order $h$ therefore contains an induced path
on at least $h-2$ vertices.

By \cref{thm:prelim-perron}, the spectral radius of a connected
graph is strictly larger than that of any proper connected subgraph.
Since every cycle has spectral radius two and $\rho(J)\le2$,
a connected component of $J$ containing a cycle must be that cycle.
The other components are trees. A vertex of degree $k$ gives a star
of spectral radius $\sqrt{k}$, so their degrees are at most four.
A tree with a vertex of degree four must be exactly $K_{1,4}$,
since this star already has spectral radius two; a path through
its center misses only two vertices.

Suppose a tree has maximum degree three and has two vertices of
degree three. Their joining path, together with the two neighbors
off the path at each end, is a double fork. The four added vertices
are distinct, because the graph is a tree. Give value two to every
vertex of the joining path and value one to its four added leaves.
The sum of neighboring values is twice the value at each vertex,
so this double fork has spectral radius two. Strict monotonicity
forces the component to equal the double fork. A path joining a leaf
at one end to a leaf at the other misses just two vertices.

A tree with exactly one vertex of degree three has three arms of
lengths $b_1\le b_2\le b_3$, measured in edges. If $b_1=1$, the path
through the other two arms misses one vertex. If $b_1\ge2$, the tree
contains the three-arm tree with all arm lengths two. Values three
at its center, two at its neighbors, and one at its leaves give a
positive eigenvector with eigenvalue two. Thus it must be the entire
component, and a path through two arms misses two vertices. Finally,
a tree of maximum degree at most two is itself a path. Every path in
a tree is induced, so a tree component of order $h$ has an induced
path on at least $h-2$ vertices.

For the construction we use the following matching calculation.
Suppose an edge form $b$ is supported on pairwise disjoint edges
$e_1,\ldots,e_q$ of an induced subgraph $X[W]$, with squared entries
$w_1,\ldots,w_q>0$. Orient the edges arbitrarily and put
$b_{uv}=\sqrt{w_i}$ on the chosen orientation of $e_i$, with the
opposite entry negative. Call two support edges \emph{conflicting}
if their four endpoints do not form a clique of $X$, equivalently
if some edge of $J$ joins their endpoints. A Pfaffian of order four can
be nonzero only on the endpoints of two support edges. Its square
there is $w_iw_j$, because exactly one of the three products in the
Pfaffian expansion is nonzero. Consequently
\begin{align}
 E_{X[W]}(b)&=\sum_i w_i,
 &T_{2,X[W]}(b)&=
 \sum_{\substack{i<j\\e_i,e_j\text{ not conflicting}}}w_iw_j,
 \label{sg:radius-matching}\\
 D_{X[W]}(b)&=\sum_i w_i^2+
 2\sum_{\substack{i<j\\e_i,e_j\text{ conflicting}}}w_iw_j.
 \label{sg:radius-conflicts}
\end{align}
In particular both $E_{X[W]}(b)$ and $D_{X[W]}(b)$ are positive
when the matching is nonempty.

Consider a component with vertex set $W$ and order $h\ge15$.
If it is a tree, choose an
induced path $v_1,\ldots,v_m$ with $m\ge h-2$ and reflect its two
ends into one another. More precisely, put
\[
 q=\left\lfloor\frac{m-1}{2}\right\rfloor,
 \qquad e_i=\{v_i,v_{m+1-i}\}\quad(1\le i\le q).
\]
These pairs are disjoint, and their indices differ by at least two,
so they are edges of $X$. For $m=2k+1$ the reflection leaves the
middle vertex unused; for $m=2k$ we omit the middle pair
$\{v_k,v_{k+1}\}$, which is an edge of $J$. Thus the retained
vertices lie in the two segments
\[
 \{v_1,\ldots,v_q\},\qquad
 \{v_{m+1-q},\ldots,v_m\}.
\]
No edge of $J$ joins these segments, and within either segment its
edges join consecutive vertices. Hence $e_i$ and $e_j$ conflict
exactly when $|i-j|=1$. All vertices off this path have zero rows
in the form, and so cannot create additional nonzero Pfaffians.

If the component is a cycle, label its vertices cyclically by
$0,1,\ldots,h-1$ and use the reflection $i\mapsto-i$ modulo $h$:
\[
 q=\left\lfloor\frac{h-2}{2}\right\rfloor,
 \qquad e_i=\{i,h-i\}\quad(1\le i\le q).
\]
For $h=2k$ this omits the fixed points $0,k$; for $h=2k+1$ it
omits the fixed point $0$ and the adjacent pair $k,k+1$.
The retained arcs $1,\ldots,q$ and $h-q,\ldots,h-1$ have no cycle
edge between them: their separation at the middle is $h-2q\ge2$,
and their separation across zero is two edges. Again each $e_i$
is an edge of $X$, and conflicts occur exactly between consecutive
indices.

In the tree case $h\le m+2\le2q+4$, and in the cycle case
$h\le2q+3$. Thus in both cases
\begin{equation}\label{sg:radius-size}
 h\le2q+4,\qquad q\ge6.
\end{equation}
Give the support edges alternating squared weights
\[
 \theta=(q+1)^{-2},\qquad
 w_i=\begin{cases}1,&i\text{ odd},\\ \theta,&i\text{ even},\end{cases}
 \qquad \ell=\lceil q/2\rceil,
 \quad \ell'=\lfloor q/2\rfloor.
\]
Thus the actual entries have magnitudes $1$ and $1/(q+1)$.
There are $q-1$ conflicting pairs, each with product $\theta$.
Writing $e=E_{X[W]}(b)$ and $\delta=D_{X[W]}(b)$ for this raw
form on the component, the matching identities give
\[
 e=\ell+\ell'\theta,
 \qquad \delta=\ell+u,
 \qquad u=\ell'\theta^2+2(q-1)\theta.
\]
Since $(q+1)^4\ge16q$ and
$(q+1)^2-8(q-1)=(q-3)^2\ge0$, we have
\[
 0\le u\le\frac{q}{2(q+1)^4}
                  +\frac{2(q-1)}{(q+1)^2}
       \le\frac1{32}+\frac14=\frac9{32}.
\]
Using $e\ge\ell$, $2\ell\ge q$, and
$\ell^2/(\ell+u)=\ell-u+u^2/(\ell+u)$, we obtain
\begin{equation}\label{sg:radius-large}
 \frac{2e^2}{\delta}
 \ge\frac{2\ell^2}{\ell+u}
 \ge2\ell-2u
 \ge q-\frac9{16}
 \ge\frac{2q+4}{3}
 \ge\frac h3.
\end{equation}
The penultimate inequality holds for $q\ge6$, since its margin
is $(16q-91)/48>0$.

For a tree of order $3\le h\le14$, the two distance-parity classes
from any root are independent sets. The larger has at least
$\lceil h/2\rceil$ vertices. Select
\[
r_T(h)=2\left\lfloor\frac{\lceil h/2\rceil}{2}\right\rfloor
\]
of them and pair these selected vertices with unit entries. They
form a clique of $X$, so there are no conflicts and the resulting
form has $e=\delta=r_T(h)/2$.

For a cycle in the same order range, the vertices
$0,2,\ldots,2(\lfloor h/2\rfloor-1)$ are independent.
Select an even number
\[
r_C(h)=2\lfloor h/4\rfloor\le\lfloor h/2\rfloor
\]
of these vertices. When this number is positive, the same unit
matching construction gives $e=\delta=r_C(h)/2$.

A direct check for $3\le h\le14$ shows that $r_T(h)\ge h/3$
throughout, while $r_C(h)\ge h/3$ except when $h=3$ or $h=7$.
Thus every nonexceptional component considered above supplies a
nonzero form with $2e^2/\delta\ge h/3$. The components still to be
treated are precisely $K_1,K_2,K_3,C_7$.

Let $W_0$ be the union of all these exceptional components, with
respective multiplicities $b_1,b_2,b_3,t_7$, and put
\[
 N_0=|W_0|=b_1+2b_2+3b_3+7t_7.
\]
For each component $K_s$, $s\in\{1,2,3\}$, take its normalized
indicator $s^{-1/2}\mathbf1_{V(K_s)}$, viewed as a column indexed
by $W_0$. For each seven-cycle, labeled cyclically by $0,\ldots,6$,
take the three coordinate columns at positions $0,2,4$. These
positions are independent in the cycle. All chosen columns have
unit norm and disjoint supports, so they are orthonormal. Their
number is
\[
 M_0=b_1+b_2+b_3+3t_7,
 \qquad 3M_0-N_0=2b_1+b_2+2t_7\ge0.
\]
Retain any $k_0=2\lfloor M_0/2\rfloor$ of the columns and assemble
them as a matrix $U_0$. Discarding at most one column in the entire
pool gives
\begin{equation}\label{sg:radius-pool}
 U_0^TU_0=I_{k_0},
 \qquad k_0\ge M_0-1\ge N_0/3-1.
\end{equation}
Put $p_0=k_0/2$ and define
\[
  b_0=U_0\Omega_{p_0}U_0^T,
\]
with $b_0=0$ if $k_0=0$. We verify the row condition of
Lemma~\ref{sg:orthonormal-construction} for $X[W_0]$.
On a component $K_s$, all rows of $U_0$ are equal, whether its
indicator column was retained or discarded. On a seven-cycle,
the only possible nonzero rows occur at the independent vertices
$0,2,4$, so every cycle edge has a zero row at one endpoint.
Thus, whenever $\{x,y\}\in E(J[W_0])$, the corresponding rows
of $U_0$ are equal or one is zero. They are therefore linearly
dependent.

The lemma shows that $b_0$ is an edge form on $X[W_0]$,
that every nonclique Pfaffian vanishes, and that
\begin{equation}\label{sg:radius-pool-energy}
 E_{X[W_0]}(b_0)=p_0,
 \qquad T_{2,X[W_0]}(b_0)=\binom{p_0}{2},
 \qquad D_{X[W_0]}(b_0)=p_0.
\end{equation}
These identities also hold for $p_0=0$; in that case this zero
block is omitted from the normalization below. Its contribution
$2p_0=k_0$ still satisfies \eqref{sg:radius-pool}.

For each nonzero raw block $b_i$ on a vertex set $W_i$, write
$e_i=E_{X[W_i]}(b_i)>0$ and
$\delta_i=D_{X[W_i]}(b_i)>0$, and set
\[
 a_i=\sqrt{e_i/\delta_i}\,b_i,
 \qquad s_i=e_i^2/\delta_i.
\]
The energy $E$ is quadratic in the entries, whereas $T_2$ and
$D$ are quartic. Hence
\begin{equation}\label{sg:radius-normalization}
 E_{X[W_i]}(a_i)=D_{X[W_i]}(a_i)=s_i>0.
\end{equation}

For each component other than $K_1,K_2,K_3$, and $C_7$, we have
$2s_i\ge h_i/3$. If
$3\le h_i\le14$, the same bound follows from the independent-set
construction, whose only exceptions among cycles are $C_3=K_3$
and $C_7$. All exceptional components have been included in $W_0$
and are handled together.

Let $a$ equal $a_i$ on each block, with all cross-block entries
zero. Every block is a union of whole components of $J$, so
distinct blocks are completely joined in $X$. In a four-vertex
principal submatrix, a block containing exactly one of the four
vertices gives a zero row. Thus a nonzero Pfaffian of order four
either lies in one block or uses two vertices from each of two
blocks. In the latter case it is, up to sign, the product of the
two internal entries; the four vertices form a clique exactly
when both internal pairs are edges of $X$. Summing gives
\begin{align*}
 E_X(a)&=\sum_iE_{X[W_i]}(a_i),\\
 T_{2,X}(a)&=\sum_iT_{2,X[W_i]}(a_i)
   +\sum_{i<j}E_{X[W_i]}(a_i)E_{X[W_j]}(a_j),\\
 D_X(a)&=\sum_iD_{X[W_i]}(a_i).
\end{align*}
The last equality follows by expanding the square in the definition
of $D_X$ and canceling the cross terms. Zero rows or an omitted
zero block do not affect these identities. By
\eqref{sg:radius-normalization}, $E_X(a)=D_X(a)=S:=\sum_i s_i$.

There is at least one nonzero block. Otherwise all components would
be exceptional and $k_0=0$, hence $M_0\le1$. A seven-cycle alone
contributes three columns, so none could occur; there would then
be at most one component, among $K_1,K_2,K_3$. This contradicts
$N\ge4$. Thus $S>0$. Finally, summing the ordinary contributions
and the exceptional contribution, including zero if necessary,
gives
\[
 2S\ge\frac{N-N_0}{3}+k_0
     \ge\frac{N-N_0}{3}+\frac{N_0}{3}-1
     =\frac N3-1.
\]
\end{proof}

\subsection{Residual classes and weighted incidence counts}

Assume throughout this subsection that $t\ge2$.
By Lemma~\ref{sg:coverage}, every vertex has two or three fat
neighbors. Put
\[
  V_i=\{v:|F(v)|=i\},\qquad n_i=|V_i|\quad(i=2,3),
  \qquad n=n_2+n_3.
\]
Since $R$ is real symmetric and positive semidefinite, it admits
a factorization $R=ZZ^T$. Let $z_v$ be the row of $Z$ indexed
by $v$. Then
\[
  R_{uv}=\langle z_u,z_v\rangle
  \qquad(u,v\in V(G)).
\]
We call the vectors $z_v$ \emph{Gram vectors} for $R$:
their pairwise inner products recover the entries of $R$.
In particular, $R_{vv}=\|z_v\|^2$.

Then
\begin{equation}\label{sg:gram-inner}
  \|z_v\|^2=3-|F(v)|,\qquad
  \langle z_u,z_v\rangle=(M_G)_{uv}-|F(u)\cap F(v)|
  \quad(u\ne v).
\end{equation}
Thus $z_v$ is a unit vector for $v\in V_2$ and is zero for
$v\in V_3$. A \emph{residual class} is an equivalence class
in $V_2$ under the relation
\[
  u\sim v\quad\Longleftrightarrow\quad z_u=\pm z_v.
\]
For each residual class $K$, choose a representative vertex
$v_K\in K$ and set $z_K:=z_{v_K}$. Define
\[
  K_+=\{u\in K:z_u=z_K\},
  \qquad
  K_-=\{u\in K:z_u=-z_K\}.
\]
For each fat vertex $f$, let $a_f$ be the number of isolated
vertices of $H_f$, let $c_f=c(H_f)$, and set $w=\max_f c_f$.

\begin{lemma}\label{sg:classes}
Distinct residual classes have orthogonal Gram directions, and
$\rank R$ is the number of classes. Within each sign side, the
fat supports are pairwise disjoint. If both sign sides are
nonempty, each has at most two vertices.

Every $H_f$ is a disjoint union of edges and isolated vertices,
and every vertex of $V_3\cap Q_f$ is isolated. Each edge $uv$ of
$H_f$ joins opposite sides of one residual class and satisfies
$F(u)\cap F(v)=\{f\}$. Consequently,
\begin{equation}\label{sg:matching-order}
  q_f+a_f=2c_f\le2w.
\end{equation}
\end{lemma}

\begin{proof}
For $u,v\in V_2$, the inner product in \eqref{sg:gram-inner} is
an integer in $[-1,1]$. Its values $\pm1$ force $z_u=\pm z_v$,
so distinct classes are orthogonal. For a residual class $K$, write $z_v=\varepsilon_vz_K$,
where $\varepsilon_v\in\{1,-1\}$ and $\|z_K\|=1$.
Then
\[
  R[K]
  =\bigl(\langle z_u,z_v\rangle\bigr)_{u,v\in K}
  =\varepsilon\varepsilon^T,
  \qquad \varepsilon=(\varepsilon_v)_{v\in K}.
\]
Every column of this matrix is a scalar multiple of the
nonzero vector $\varepsilon$, so its rank is exactly one.
Distinct classes are orthogonal, and the Gram vectors on
$V_3$ are zero. Thus, after reordering the vertices, $R$ is
block diagonal with one rank-one block for each residual
class and a zero block on $V_3$. Consequently, $\rank R$
equals the number of residual classes.

Distinct vertices with the same Gram vector are adjacent and have
disjoint fat supports. Vertices with opposite Gram vectors have
intersecting supports. Suppose both $K_+$ and $K_-$ are nonempty, and fix $u\in K_+$.
Every $v\in K_-$ has a Gram vector opposite to that of $u$,
so $F(v)\cap F(u)\ne\varnothing$. On the other hand, the
sets $F(v)$, $v\in K_-$, are pairwise disjoint, because these
vertices have the same Gram vector. Hence the intersections
$F(v)\cap F(u)$ are nonempty, pairwise disjoint subsets of
the two-element set $F(u)$. Therefore
\[
  |K_-|
  \le \sum_{v\in K_-}|F(v)\cap F(u)|
  \le |F(u)|=2.
\]
Interchanging $K_+$ and $K_-$ gives $|K_+|\le2$ as well.

If $uv\in E(H_f)$, then
$\langle z_u,z_v\rangle=-|F(u)\cap F(v)|\le-1$.
The Gram norms force $u,v\in V_2$, $z_u=-z_v$, and
$F(u)\cap F(v)=\{f\}$. Two neighbors of a vertex in $H_f$
would therefore have the same Gram vector and disjoint fat
supports, although both contain $f$. Hence $H_f$ has maximum
degree at most one. Its edges and isolated vertices give
$q_f+a_f=2c_f$, proving \eqref{sg:matching-order}.
\end{proof}

We now prove a technical lemma based on combinatorial counting.

Call an incidence $(v,f)$ with $f\in F(v)$ \emph{unpaired} if
$v$ is isolated in $H_f$, and \emph{paired} otherwise. For
$v\in V_2$, let $u_v$ be its number of unpaired incidences, and put
\[
  U_2=\sum_{v\in V_2}u_v,\qquad
  P_2=\sum_{v\in V_2}
      \sum_{\substack{f\in F(v)\\(v,f)\text{ unpaired}}}q_f.
\]
Thus $0\le U_2\le2n_2$. For a residual class $K$, let $U_K,P_K$
denote the same sums restricted to $v\in K$. Finally, set
\[
  s_*=d-2w-B_*,\qquad \mathcal B=4w-d+B_*-4.
\]

\begin{starredlemma}\label{sg:counts}The incidence counts satisfy\begin{align}  q_f&\ge s_*+2\quad\text{if }Q_f\cap V_2\ne\varnothing,    \label{sg:size}\\  (s_*+2)U_2\le P_2&\le\mathcal Bn+U_2-(s_*+2B_*+2)n_3,    \label{sg:weighted}\\  P_2&\le(\mathcal B+2)n\quad\text{if }w\le d/3.    \label{sg:P-bound}\end{align}\end{starredlemma}

\begin{proof}
Let $v\in V_2$. By the definition of $V_2$, the vertex $v$ has
exactly two fat neighbors, which we denote by $f$ and $g$.
Thus $F(v)=\{f,g\}$, and the two incidences at $v$ are
$(v,f)$ and $(v,g)$. Since $u_v$ counts the unpaired incidences
at $v$, exactly $2-u_v$ of these two incidences are paired.
We first show that these paired incidences correspond bijectively
to the vertices that share a fat neighbor with $v$ but are
nonadjacent to $v$ in $G$.

Suppose that $(v,f)$ is paired. Then $v$ is not isolated in
$H_f$. By Lemma~\ref{sg:classes}, it has a unique neighbor
$y_f$ in $H_f$, and
\[
  F(v)\cap F(y_f)=\{f\}.
\]
Since $H_f=\overline{G[Q_f]}$, the vertices $v$ and $y_f$
are nonadjacent in $G$. The displayed identity also gives
$y_f\in Q_f\setminus Q_g$. Similarly, if
$(v,g)$ is paired, it determines a unique vertex
$y_g\in Q_g\setminus Q_f$ that is nonadjacent to $v$ in $G$.
If both incidences are paired, their corresponding vertices
are therefore distinct.

Conversely, let $y\ne v$ be nonadjacent to $v$ in $G$ and
suppose that $F(y)\cap F(v)\ne\varnothing$. Choose
$h\in F(y)\cap F(v)$; necessarily $h\in\{f,g\}$.
Then $v,y\in Q_h$, and their nonadjacency in $G$ means that
$vy\in E(H_h)$. Thus $(v,h)$ is paired and $y$ is its unique
corresponding vertex. This proves the claimed bijection.
Consequently, exactly $2-u_v$ vertices of
$(Q_f\cup Q_g)\setminus\{v\}$ are nonadjacent to $v$ in $G$,
and each belongs to exactly one of $Q_f$ and $Q_g$.

By inclusion--exclusion and Lemma~\ref{sg:coverage},
\[
  \begin{aligned}
    q_f+q_g-2
    &=|Q_f\setminus\{v\}|+|Q_g\setminus\{v\}|\\
    &=(2-u_v)+d-|N_u(v)|
      +|(Q_f\cap Q_g)\setminus\{v\}|\\
    &\ge d+2-B_*-u_v.
  \end{aligned}
\]
Adding two to both sides yields
\begin{equation}\label{sg:incidence-degree-lower}
  q_f+q_g\ge d+4-B_*-u_v.
\end{equation}
Using $u_v\le2$ and $q_g\le2w$, the latter following from
\eqref{sg:matching-order} in Lemma~\ref{sg:classes}, we obtain
\[
  \begin{aligned}
    q_f
    &\ge d+4-B_*-u_v-q_g\\
    &\ge d+4-B_*-2-2w\\
    &=s_*+2.
  \end{aligned}
\]
For any fat vertex $f$ with $Q_f\cap V_2\ne\varnothing$,
we may choose $v$ in that intersection and apply the same
argument. This proves \eqref{sg:size}.
In particular, every unpaired incidence $(v,h)$ with
$v\in V_2$ has $q_h\ge s_*+2$. Summing this bound gives
\[
  \begin{aligned}
    P_2
    &=\sum_{v\in V_2}
      \sum_{\substack{h\in F(v)\\(v,h)\text{ unpaired}}}q_h\\
    &\ge (s_*+2)\sum_{v\in V_2}u_v
     =(s_*+2)U_2.
  \end{aligned}
\]

Next let $v\in V_3$. Since $z_v=0$, the Gram identity
\eqref{sg:gram-inner} implies, for every $u\ne v$,
\[
  \begin{aligned}
    0=\langle z_u,z_v\rangle
    &=(M_G)_{uv}-|F(u)\cap F(v)|,\\
    \text{hence}\qquad |F(u)\cap F(v)|
    &=(M_G)_{uv}\in\{0,1\}.
  \end{aligned}
\]
Thus $u$ is adjacent to $v$ if and only if they share a fat
neighbor, and in that case the shared fat neighbor is unique.
Every neighbor of $v$ therefore lies in exactly one of the
sets $Q_f\setminus\{v\}$ with $f\in F(v)$, and every vertex
in any of these sets is a neighbor of $v$. These sets form
a disjoint partition of $N_G(v)$. Since $|F(v)|=3$ and
$|N_G(v)|=d$, it follows that
\[
  d=\sum_{f\in F(v)}(q_f-1)
   =\sum_{f\in F(v)}q_f-3,
  \qquad
  \sum_{f\in F(v)}q_f=d+3.
\]
Moreover, by Lemma~\ref{sg:classes}, the vertex $v\in V_3$
is isolated in $H_f$ for every $f\in F(v)$. Hence all three
incidences at $v$ are unpaired.

We now count the weighted incidences over the whole graph.
Below, $\sum_f$ denotes summation over $V_f$.
For a fixed fat vertex $f$, there are $q_f$ vertices in $Q_f$,
each contributing weight $q_f$. Thus
\[
  \sum_f q_f^2
  =\sum_f\sum_{v\in Q_f}q_f
  =\sum_{v\in V(G)}\sum_{f\in F(v)}q_f.
\]
The last equality counts the same incidences in the opposite
order: first choosing the slim vertex and then its fat neighbor.
For $v\in V_2$, \eqref{sg:incidence-degree-lower} bounds its
inner sum below by $d+4-B_*-u_v$; for $v\in V_3$, the inner
sum is exactly $d+3$. Since $V(G)=V_2\sqcup V_3$, we obtain
\[
  \begin{aligned}
    \sum_f q_f^2
    &=\sum_{v\in V_2}\sum_{f\in F(v)}q_f
      +\sum_{v\in V_3}\sum_{f\in F(v)}q_f\\
    &\ge\sum_{v\in V_2}(d+4-B_*-u_v)
      +\sum_{v\in V_3}(d+3)\\
    &=(d+4-B_*)n_2-U_2+(d+3)n_3.
  \end{aligned}
\]

Similarly, $a_f$ counts precisely the vertices $v\in Q_f$
for which $(v,f)$ is unpaired. Giving each such incidence
weight $q_f$ and again reversing the order of summation gives
\[
  \begin{aligned}
    \sum_f q_fa_f
    &=\sum_f
      \sum_{\substack{v\in Q_f\\(v,f)\text{ unpaired}}}q_f\\
    &=\sum_{v\in V(G)}
      \sum_{\substack{f\in F(v)\\(v,f)\text{ unpaired}}}q_f\\
    &=P_2+(d+3)n_3.
  \end{aligned}
\]
Here the contribution from $V_2$ is $P_2$ by definition.
For each $v\in V_3$, all its incidences are unpaired and their
weights sum to $d+3$, which accounts for the second term.
Adding the preceding two estimates yields
\[
  \sum_f q_f(q_f+a_f)
  \ge(d+4-B_*)n_2-U_2+P_2+2(d+3)n_3.
\]

For an upper bound on the same sum, first count all incidences
without weights. Each vertex of $V_2$ has two fat neighbors
and each vertex of $V_3$ has three, so
\[
  \sum_f q_f
  =\sum_{v\in V(G)}|F(v)|
  =2n_2+3n_3.
\]
By \eqref{sg:matching-order} in Lemma~\ref{sg:classes}, we have
$q_f+a_f=2c_f\le2w$ for every $f$. Multiplying by $q_f\ge0$
and summing gives
\[
  \sum_f q_f(q_f+a_f)
  \le2w\sum_f q_f
  =2w(2n_2+3n_3).
\]
Combining the lower and upper bounds, we have
\[
  (d+4-B_*)n_2-U_2+P_2+2(d+3)n_3
  \le2w(2n_2+3n_3).
\]
Isolating $P_2$ and collecting the coefficients of $n_2$ and
$n_3$ gives
\[
  \begin{aligned}
    P_2
    &\le 2w(2n_2+3n_3)-(d+4-B_*)n_2
         +U_2-2(d+3)n_3\\
    &=(4w-d+B_*-4)n_2+(6w-2d-6)n_3+U_2\\
    &=\mathcal B(n-n_3)+(6w-2d-6)n_3+U_2\\
    &=\mathcal Bn+U_2-(d-2w+B_*+2)n_3\\
    &=\mathcal Bn+U_2-(s_*+2B_*+2)n_3.
  \end{aligned}
\]
In the last three lines we used $n_2=n-n_3$,
$\mathcal B=4w-d+B_*-4$, and $s_*=d-2w-B_*$.
Together with the lower bound $P_2\ge(s_*+2)U_2$ proved
above, this establishes \eqref{sg:weighted}.

Finally, suppose $w\le d/3$. Then
\[
  s_*+2B_*+2
  =d-2w+B_*+2
  \ge\frac d3+B_*+2>0.
\]
Since $n_3\ge0$, the term $-(s_*+2B_*+2)n_3$ is
nonpositive and may be discarded when taking an upper bound.
Also, each vertex in $V_2$ contributes at most two unpaired
incidences, so $U_2\le2n_2\le2n$. Consequently,
\[
  \begin{aligned}
    P_2
    &\le\mathcal Bn+U_2-(s_*+2B_*+2)n_3\\
    &\le\mathcal Bn+U_2\\
    &\le\mathcal Bn+2n
     =(\mathcal B+2)n,
  \end{aligned}
\]
which proves \eqref{sg:P-bound}.
\end{proof}

\subsection{Rank estimates and completion by a local edge form}

For the following rank estimates, retain the assumption $t\ge2$.
We use the bounds from Lemma~\ref{sg:counts} to estimate
the rank of $M_G+3I=R+CC^T$ in both directions.

\begin{starredlemma}\label{sg:rank-bounds}Assume $w\le d/3$. Then\begin{equation}\label{sg:rank-forward}  \eta\le\frac{3w}{d+2}-\frac12       +\frac{3B_*}{4(d+2)}+\frac3{s_*+2},\end{equation}and\begin{equation}\label{sg:rank-reverse}  \eta\ge  \frac{2d-6w-2B_*+5}{4(s_*+1)}.\end{equation}In particular, if $0\le\eta\le1/4$, then\begin{equation}\label{sg:rank-invert}  w\ge  \frac{(1-2\eta)(d-B_*)+(5-4\eta)/2}{3-4\eta}.\end{equation}\end{starredlemma}

\noindent\emph{Readers may proceed directly to the \hyperref[sg:proof-theorem-b]{proof of Theorem~\ref*{sg:theorem}}.}

\begin{proof}
Our standing degree assumption and $w\le d/3$ give
\begin{equation}\label{sg:rank-positive}
  s_*=d-2w-B_*\ge\frac d3-B_*\ge\frac d6>0.
\end{equation}
Thus both denominators $s_*+1$ and $s_*+2$ are positive.

We first prove the forward bound.

Fix a residual class $K$ with $r=|K_+|$, $s=|K_-|$, and $k=r+s$.
For each $u\in K$, write $z_u=\varepsilon_u z_K$, where
$\varepsilon_u=1$ on $K_+$ and $\varepsilon_u=-1$ on $K_-$.
For fixed $v\in K$, Lemma~\ref{sg:classes} and the vanishing
of the Gram vectors on $V_3$ give
\[
  R_{vu}=\langle z_v,z_u\rangle
  =\begin{cases}
    \varepsilon_v\varepsilon_u,&u\in K,\\
    0,&u\notin K.
  \end{cases}
\]
Also,
\[
  (CC^T)_{vu}=\sum_{f\in V_f}C_{vf}C_{uf},
  \qquad
  \sum_{u\in V(G)}C_{uf}=q_f.
\]
Since $G$ is $d$-regular, summing the entries in row $v$ of
$M_G+3I=R+CC^T$ yields
\[
  \begin{aligned}
    d+3
    &=\sum_{u\in V(G)}(M_G+3I)_{vu}\\
    &=\sum_{u\in V(G)}R_{vu}
      +\sum_{u\in V(G)}(CC^T)_{vu}\\
    &=\varepsilon_v\sum_{u\in K}\varepsilon_u
      +\sum_{f\in V_f}C_{vf}\sum_{u\in V(G)}C_{uf}\\
    &=\varepsilon_v(r-s)+\sum_{f\in F(v)}q_f.
  \end{aligned}
\]
Hence
\[
  \sum_{f\in F(v)}q_f=d+3-\varepsilon_v(r-s)
  \qquad(v\in K).
\]

 Summing over $v\in K$ and using
$\sum_{v\in K}\varepsilon_v=r-s$, we obtain the total incidence
weight
\[
  \begin{aligned}
    \sum_{v\in K}\sum_{f\in F(v)}q_f
    &=\sum_{v\in K}\bigl(d+3-\varepsilon_v(r-s)\bigr)\\
    &=k(d+3)-(r-s)\sum_{v\in K}\varepsilon_v\\
    &=k(d+3)-(r-s)^2.
  \end{aligned}
\]

Now consider a paired incidence $(v,f)$ with $v\in K$.
By Lemma~\ref{sg:classes}, the vertex $v$ has a unique
neighbor $u$ in $H_f$, and $u$ lies on the opposite sign side
of the same residual class $K$. The same lemma gives
$F(u)\cap F(v)=\{f\}$. Thus the paired incidences with slim
vertices in $K$ split into disjoint groups
\[
  \{(u,f),(v,f)\},
  \qquad u\in K_+,\quad v\in K_-.
\]
Each group determines a vertex pair in $K_+\times K_-$.
Since $F(u)\cap F(v)=\{f\}$, the same vertex pair cannot
occur in groups associated with two different fat vertices.
There are therefore at most $|K_+||K_-|=rs$ groups.
By \eqref{sg:matching-order} in Lemma~\ref{sg:classes}, each
such group has weight $2q_f\le4w$. Hence
\[
  \sum_{v\in K}
    \sum_{\substack{f\in F(v)\\(v,f)\text{ paired}}}q_f
  \le4wrs.
\]

Since $P_K$ counts only the unpaired incidence weight, subtracting
this bound from the total weight gives
\[
  \begin{aligned}
    P_K
    &=\sum_{v\in K}\sum_{f\in F(v)}q_f
      -\sum_{v\in K}
        \sum_{\substack{f\in F(v)\\(v,f)\text{ paired}}}q_f\\
    &\ge k(d+3)-(r-s)^2-4wrs.
  \end{aligned}
\]
We claim $3P_K\ge(4-k)(d+2)$. This is immediate for $k\ge4$.
For $k\le3$, either $rs=0$, or $rs=k-1$ and $(r-s)^2=(k-2)^2$.
Using $w\le d/3$, respectively, gives
\[
  3P_K-(4-k)(d+2)\ge
  \begin{cases}
    (k-1)(4d+8-3k),&rs=0,\\
    (k-1)(20-3k),&rs>0.
  \end{cases}
\]
Both expressions are nonnegative. Together with the case $k\ge4$,
this proves that every residual class $K$ satisfies
\[
  3P_K\ge(4-|K|)(d+2).
\]
Let $m_R$ be the number of residual classes. By
Lemma~\ref{sg:classes}, $m_R=\rank R$.
The residual classes partition $V_2$, and $P_K$ is the
contribution of $K$ to $P_2$. Hence
\[
  \sum_K |K|=n_2,
  \qquad
  \sum_K P_K=P_2,
\]
where both sums run over all residual classes. Summing the
preceding inequality over these $m_R$ classes gives
\[
  \begin{aligned}
    3P_2
    &=\sum_K 3P_K\\
    &\ge(d+2)\sum_K(4-|K|)\\
    &=(d+2)(4m_R-n_2).
  \end{aligned}
\]
Rearranging and using $m_R=\rank R$, we obtain
\begin{equation}\label{sg:rank-class-upper}
  \rank R\le\frac{n_2}{4}+\frac{3P_2}{4(d+2)}.
\end{equation}

To bound the contribution of $CC^T$ to the rank, define
\[
  \begin{aligned}
    F_L&:=\bigl|\{f\in V_f:Q_f\cap V_2\ne\varnothing\}\bigr|,\\
    F_S&:=\bigl|\{f\in V_f:Q_f\cap V_2=\varnothing\}\bigr|.
  \end{aligned}
\]
These two classes partition the fat vertices, so
$|V_f|=F_L+F_S$.

By \eqref{sg:size} in Lemma~\ref{sg:counts}, every fat
neighborhood meeting $V_2$ has order at least $s_*+2$.
Counting the incidences belonging to these neighborhoods and
then bounding them by all incidences gives
\[
  \begin{aligned}
    (s_*+2)F_L
    &\le\sum_{\substack{f\in V_f\\Q_f\cap V_2\ne\varnothing}}q_f\\
    &\le\sum_{f\in V_f}q_f
     =\sum_{v\in V(G)}|F(v)|\\
    &=2n_2+3n_3\le3n.
  \end{aligned}
\]
Since $s_*+2>0$ by \eqref{sg:rank-positive}, it follows that
\[
  F_L\le\frac{2n_2+3n_3}{s_*+2}
      \le\frac{3n}{s_*+2}.
\]

Every fat neighborhood contains a clique of order at least
$r_*$ by the associated-graph construction and
\eqref{sg:clique-covered}, so $q_f\ge r_*$.
If $Q_f\cap V_2=\varnothing$, then $Q_f\subseteq V_3$.
Thus all incidences counted by these neighborhoods have their
slim vertices in $V_3$. Since each vertex of $V_3$ has three
fat neighbors,
\[
  r_*F_S
  \le\sum_{\substack{f\in V_f\\Q_f\cap V_2=\varnothing}}q_f
  \le\sum_{v\in V_3}|F(v)|
  =3n_3.
\]
Consequently, using $r_*\ge169$, we obtain
\[
  F_S\le\frac{3n_3}{r_*}\le\frac{n_3}{4}.
\]

The matrix $C$ has one column for each fat vertex, hence
$F_L+F_S$ columns. The column space of $CC^T$ is contained
in that of $C$, and the column space of a matrix sum is
contained in the sum of the two column spaces. Therefore
\[
  \begin{aligned}
    \rank(M_G+3I)
    &=\rank(R+CC^T)\\
    &\le\rank R+\rank(CC^T)\\
    &\le\rank R+\rank C\\
    &\le\rank R+F_L+F_S.
  \end{aligned}
\]
Combining the bounds for $F_L,F_S$ with
\eqref{sg:rank-class-upper} gives
\[
  \begin{aligned}
    \rank(M_G+3I)
    &\le\frac{n_2}{4}+\frac{3P_2}{4(d+2)}
      +\frac{3n}{s_*+2}+\frac{n_3}{4}\\
    &=\frac n4+\frac{3P_2}{4(d+2)}+\frac{3n}{s_*+2}.
  \end{aligned}
\]
Finally, divide by $n$ and use \eqref{sg:P-bound} to obtain
\[
  \begin{aligned}
    \eta
    &\le\frac14+\frac{3P_2}{4n(d+2)}+\frac3{s_*+2}\\
    &\le\frac14+\frac{3(\mathcal B+2)}{4(d+2)}
      +\frac3{s_*+2}\\
    &=\frac{3w}{d+2}-\frac12
      +\frac{3B_*}{4(d+2)}+\frac3{s_*+2}.
  \end{aligned}
\]

This proves \eqref{sg:rank-forward}.

For the reverse bound, fix a residual class $K$ and recall that
$U_K=\sum_{v\in K}u_v$ counts its unpaired incidences.
We show that
\[
  \frac{|K|}{4}-\frac{U_K}{8}\le1
\]
by considering the two possible configurations of its sign sides.

If both $K_+$ and $K_-$ are nonempty, Lemma~\ref{sg:classes}
gives $|K_+|\le2$ and $|K_-|\le2$, hence $|K|\le4$.
Since $U_K\ge0$, we obtain
\[
  \frac{|K|}{4}-\frac{U_K}{8}
  \le\frac{|K|}{4}\le1.
\]

If only one sign side is nonempty, all Gram vectors in $K$
are equal. In this case every incidence $(v,f)$ with
$v\in K$ and $f\in F(v)$ is unpaired. Indeed, if $(v,f)$
were paired, then $v$ would have a neighbor $u$ in $H_f$.
By Lemma~\ref{sg:classes}, the vertices $u,v$ would belong
to opposite sign sides of the same residual class, contradicting
the assumption on $K$. Each $v\in K\subseteq V_2$ has exactly
two fat neighbors, so $u_v=2$ and
\[
  U_K=\sum_{v\in K}u_v=2|K|,
  \qquad
  \frac{|K|}{4}-\frac{U_K}{8}
  =\frac{|K|}{4}-\frac{2|K|}{8}=0\le1.
\]

Write the residual classes as $K_1,\ldots,K_{m_R}$.
By Lemma~\ref{sg:classes}, $m_R=\rank R$.
Since these classes partition $V_2$, and $U_K$ is the
contribution of $K$ to $U_2$, we have
\[
  \sum_{i=1}^{m_R}|K_i|=n_2,
  \qquad
  \sum_{i=1}^{m_R}U_{K_i}=U_2.
\]
Summing the inequality just proved over all residual classes gives
\[
  \begin{aligned}
    \rank R
    &=m_R=\sum_{i=1}^{m_R}1\\
    &\ge\sum_{i=1}^{m_R}
      \left(\frac{|K_i|}{4}-\frac{U_{K_i}}{8}\right)\\
    &=\frac{n_2}{4}-\frac{U_2}{8}.
  \end{aligned}
\]
Since $R$ and $CC^T$ are positive semidefinite,
$\ker(R+CC^T)=\ker R\cap\ker C^T$, so
$\rank(M_G+3I)\ge\rank R$. Therefore
\begin{equation}\label{sg:rank-unpaired-lower}
  \eta\ge\frac14-\frac{U_2+2n_3}{8n}.
\end{equation}
By \eqref{sg:weighted},
$(s_*+1)U_2+(s_*+2B_*+2)n_3\le\mathcal Bn$.
Hence
\[
  (s_*+1)(U_2+2n_3)
  \le2\bigl((s_*+1)U_2+(s_*+2B_*+2)n_3\bigr)
  \le2\mathcal Bn.
\]
Inserting this in \eqref{sg:rank-unpaired-lower} gives
\[
  \eta\ge\frac14-\frac{\mathcal B}{4(s_*+1)}
  =\frac{2d-6w-2B_*+5}{4(s_*+1)},
\]
which is \eqref{sg:rank-reverse}. Rearrangement yields
\[
  (3-4\eta)w
  \ge(1-2\eta)(d-B_*)+\frac{5-4\eta}{2}.
\]
For $0\le\eta\le1/4$, division by $3-4\eta>0$ proves
\eqref{sg:rank-invert}.
\end{proof}

\phantomsection\label{sg:proof-theorem-b}
\begin{proof}[\textnormal{\textbf{Proof of Theorem~\ref*{sg:theorem}}}]\leavevmode\par\nobreak
By Lemma~\ref{sg:coverage}, $t\in\{1,2,3\}$.
We construct a nonzero form $a_Q$ on some
$G[Q_f]$, with $D_{G[Q_f]}(a_Q)>0$, such that
\begin{equation}\label{sg:local-target}
  \frac{2E_{G[Q_f]}(a_Q)^2}{D_{G[Q_f]}(a_Q)}
  \ge d\,g(\eta)-C_*.
\end{equation}

If $t=1$, choose a vertex with unique fat neighbor $f$.
By \eqref{sg:fat-average} and Lemma~\ref{sg:coverage},
$q_f\ge d-B_*+1\ge6$ and $\rho(H_f)\le2$.
Lemma~\ref{sg:radius} supplies the required nonzero form, since
\[
  \frac{2E_{G[Q_f]}(a_Q)^2}{D_{G[Q_f]}(a_Q)}
  \ge\frac{q_f}{3}-1
  \ge\frac d3-\frac{B_*}{3}-\frac23
  \ge d\,g(\eta)-C_*.
\]

If $t\ge2$, choose $f$ with $c_f=w$ and set
$p_0=\lfloor w/2\rfloor$. Since
$w\ge q_f/2\ge r_*/2$, we have $p_0>0$.
The form from Lemma~\ref{sg:kernel} satisfies
\begin{equation}\label{sg:completion-kernel-value}
  E_{G[Q_f]}(a_Q)=D_{G[Q_f]}(a_Q)=p_0,
  \qquad
  \frac{2E_{G[Q_f]}(a_Q)^2}{D_{G[Q_f]}(a_Q)}
  =h:=2p_0\ge w-1.
\end{equation}
If $w\ge(d-B_*+1)/3$, the preceding bound by
$d/3-B_*/3-2/3$ again proves \eqref{sg:local-target}.
It remains to consider
\begin{equation}\label{sg:completion-small-w}
  w<\frac{d-B_*+1}{3}\le\frac d3.
\end{equation}
Lemma~\ref{sg:rank-bounds} therefore applies.

For $0\le\eta\le1/4$, \eqref{sg:rank-invert} gives
\[
  \begin{aligned}
    h&\ge w-1\\
     &\ge d\,g(\eta)-B_*g(\eta)
          +\frac{4\eta-1}{2(3-4\eta)}\\
     &\ge d\,g(\eta)-\frac{B_*}{3}-\frac16
      \ge d\,g(\eta)-C_*.
  \end{aligned}
\]
For $1/4\le\eta\le1$, \eqref{sg:rank-forward} gives
\[
  \begin{aligned}
    h&\ge\frac{(2\eta+1)(d+2)}6-\frac{B_*}{4}
          -\frac{d+2}{s_*+2}-1\\
     &\ge d\,g(\eta)-\frac{B_*}{4}-7
      \ge d\,g(\eta)-C_*.
  \end{aligned}
\]
Here we used $g(\eta)\le(2\eta+1)/6$ on this interval,
$(d+2)/(s_*+2)\le6$ by \eqref{sg:rank-positive}, and
$B_*/4+7\le2B_*+6$.
Thus \eqref{sg:local-target} holds in every case.

Extend $a_Q$ by zero to an edge form $a$ on $G$. A clique
meeting $V(G)\setminus Q_f$ has a zero row in its restricted
matrix, so
\begin{equation}\label{sg:completion-extension}
  T_{j,G}(a)=T_{j,G[Q_f]}(a_Q)\quad(j\ge1),
  \qquad D_G(a)=D_{G[Q_f]}(a_Q)>0.
\end{equation}
This proves \eqref{sg:certificate-target}. Finally,
\eqref{eq:certificate-introduction} gives
\[
  d\,g(\eta)-C_*
  \le\frac{2E_G(a)^2}{D_G(a)}
  \le2\left\lfloor\frac{\omega(G)}2\right\rfloor
  \le\omega(G).
\]
All counts are taken over the whole graph, and the Hoffman
construction applies componentwise, so the proof also covers
disconnected regular graphs.
\end{proof}

\subsection{Optimality}
\label{sg:comparison-optimality}
Since $g(\eta)\ge1/4$, Theorem~\ref{sg:theorem} gives the uniform bound
\[
\omega(G)\ge \frac d4-C_*.
\]
The coefficient $1/4$ is optimal by the construction of
Yang and Koolen~\cite[Section~1]{YangKoolen2025}:
\begin{proposition}
\label{sg:optimality-construction}
For arbitrarily large $v\ge13$, the construction of Yang and Koolen
from a Steiner system $S(2,4,v)$ yields a connected $d_v$-regular
graph $G_v$ such that, with $r=(v-1)/3$,
\[
d_v=4r-3,\qquad
\lambda_{\min}(M_{G_v})\ge-3,\qquad
\omega(G_v)=r=\frac{d_v+3}{4}\sim\frac{d_v}{4}
\quad(v\to\infty).
\]
\end{proposition}
\begin{proof}
See~\cite[Section~1]{YangKoolen2025} for the construction,
regularity and spectral bound. Its vertices are the four cyclic
pairs in each Steiner block; within a block only opposite pairs
are adjacent, and between blocks adjacency means intersection.
A clique containing an opposite pair cannot have a third vertex,
since distinct Steiner blocks intersect in at most one point.
Any larger clique therefore uses distinct blocks and has pairwise
intersecting two-element labels. Such a family either has a common
point or has at most three members. In the first case its size is
at most the number $r$ of blocks through that point; hence every
clique has size at most $r\ge4$.
Conversely, choosing one cyclic pair through a fixed point in each
incident block gives a clique of size $r$. Thus $\omega(G_v)=r$.
\end{proof}
Thus no bound $\omega(G)\ge c\,d-C$ with a graph-independent
constant $C$ can hold for this class with $c>1/4$.

\section{Further questions}
\label{sec:questions}

The spectral application in this paper uses only the case $j=1$,
$k=2$ of the helicity energy Maclaurin inequality. Can the higher-order
inequalities yield stronger clique certificates or further
structural information about the graph? Although the constants
in the general inequalities are sharp, the higher-order energies
may capture information beyond the two energies $E_G(a)$ and
$T_{2,G}(a)$ alone. A natural problem is to construct edge forms whose
higher-order energies can be estimated from spectral or
combinatorial data, and to determine whether these estimates
lead to stronger conclusions under additional hypotheses.

A second question concerns the relationship between graph
helicity and the cohomology of the clique complex. Under what
conditions can the local helicities be assembled into global
functionals with a cohomological interpretation? A general
clique complex need not carry a fundamental class, so this
question also requires identifying suitable cycles on which
to evaluate helicity. It would be interesting to understand
how these functionals interact with the cohomological structure
of the complex and what information they provide beyond its
cohomology groups.

More broadly, can helicity be developed into a systematic
method for combinatorial problems, analogous to the entropy
method? The present work gives one instance in which local
helicity identities lead to algebraic inequalities and then
to a graph-theoretic bound. Developing this approach further
would require general principles for choosing edge forms,
estimating their helicity energies, and translating these
estimates into combinatorial conclusions. Applications beyond
clique bounds would help clarify the scope of such a method.

\section*{Acknowledgements and statement on AI assistance}

The author thanks Kaibo Hu for introducing this project and for many valuable suggestions on its research directions, and Puchun Zhou for his generous help with the writing of this paper.

The author originated the helicity energy Maclaurin inequality and the Pfaffian realization of discrete helicity that motivated this paper. The complete proof of the rank-dependent clique-number result, Theorem~\ref{sg:theorem}, was produced entirely by GPT Astra ultra. The author has reviewed all proofs. In preparing the present manuscript, GPT Astra ultra also assisted with English exposition and organization, using those reviewed proofs.

\bibliographystyle{plain}
\bibliography{references}

@phdthesis{Hirani2003,
  author  = {Hirani, Anil N.},
  title   = {Discrete Exterior Calculus},
  school  = {California Institute of Technology},
  address = {Pasadena, California},
  year    = {2003},
  doi     = {10.7907/ZHY8-V329},
  url     = {https://resolver.caltech.edu/CaltechETD:etd-05202003-095403}
}

@misc{DesbrunHiraniLeokMarsden2005,
  author        = {Desbrun, Mathieu and Hirani, Anil N. and Leok, Melvin and Marsden, Jerrold E.},
  title         = {Discrete Exterior Calculus},
  year          = {2005},
  eprint        = {math/0508341},
  archivePrefix = {arXiv},
  primaryClass  = {math.DG},
  note          = {arXiv:math/0508341},
  doi           = {10.48550/arXiv.math/0508341},
  url           = {https://arxiv.org/abs/math/0508341}
}

@article{Kotiuga1989,
  author  = {Kotiuga, P. R.},
  title   = {Helicity functionals and metric invariance in three dimensions},
  journal = {IEEE Transactions on Magnetics},
  volume  = {25},
  number  = {4},
  pages   = {2813--2815},
  year    = {1989},
  doi     = {10.1109/20.34293},
  url     = {https://doi.org/10.1109/20.34293}
}

@article{Kotiuga1990,
  author  = {Kotiuga, P. R.},
  title   = {Analysis of finite-element matrices arising from discretizations of helicity functionals},
  journal = {Journal of Applied Physics},
  volume  = {67},
  number  = {9},
  pages   = {5815--5817},
  year    = {1990},
  doi     = {10.1063/1.345973},
  url     = {https://doi.org/10.1063/1.345973}
}

@book{Whitney1957,
  author    = {Whitney, Hassler},
  title     = {Geometric Integration Theory},
  publisher = {Princeton University Press},
  address   = {Princeton, New Jersey},
  year      = {1957}
}

@book{HLP1952,
  author    = {Hardy, G. H. and Littlewood, J. E. and P{\'o}lya, G.},
  title     = {Inequalities},
  edition   = {2},
  publisher = {Cambridge University Press},
  address   = {Cambridge},
  year      = {1952}
}

@article{MotzkinStraus1965,
author  = {Motzkin, T. S. and Straus, E. G.},
title   = {Maxima for graphs and a new proof of a theorem of {Tur{\'a}n}},
journal = {Canadian Journal of Mathematics},
volume  = {17},
pages   = {533--540},
year    = {1965},
doi     = {10.4153/CJM-1965-053-6},
url     = {https://doi.org/10.4153/CJM-1965-053-6}
}

@article{FisherRyan1992,
  author  = {Fisher, David C. and Ryan, Jennifer},
  title   = {Bounds on the number of complete subgraphs},
  journal = {Discrete Mathematics},
  volume  = {103},
  number  = {3},
  pages   = {313--320},
  year    = {1992},
  doi     = {10.1016/0012-365X(92)90323-8},
  url     = {https://doi.org/10.1016/0012-365X(92)90323-8}
}

@misc{Nikiforov2006,
  author        = {Nikiforov, Vladimir},
  title         = {An extension of {Maclaurin}'s inequalities},
  year          = {2006},
  eprint        = {math/0608199},
  archivePrefix = {arXiv},
  primaryClass  = {math.CO},
  note          = {arXiv:math/0608199},
  url           = {https://arxiv.org/abs/math/0608199}
}

@article{AragaoSouza2024,
  author  = {Arag{\~a}o, Lucas and Souza, Victor},
  title   = {Localised graph {Maclaurin} inequalities},
  journal = {Annals of Combinatorics},
  volume  = {28},
  pages   = {1021--1033},
  year    = {2024},
  doi     = {10.1007/s00026-023-00672-0},
  url     = {https://doi.org/10.1007/s00026-023-00672-0}
}

@article{KYY,
  author        = {Koolen, Jack H. and Yang, Jae Young and Yang, Qianqian},
  title         = {On graphs with smallest eigenvalue at least {$-3$} and their lattices},
  journal       = {Advances in Mathematics},
  volume        = {338},
  pages         = {847--864},
  year          = {2018},
  doi           = {10.1016/j.aim.2018.09.004},
  eprint        = {1804.00369},
  archivePrefix = {arXiv},
  primaryClass  = {math.CO},
  url           = {https://doi.org/10.1016/j.aim.2018.09.004}
}

@book{BH2012,
  author    = {Brouwer, Andries E. and Haemers, Willem H.},
  title     = {Spectra of Graphs},
  series    = {Universitext},
  publisher = {Springer},
  address   = {New York},
  year      = {2012},
  doi       = {10.1007/978-1-4614-1939-6},
  isbn      = {978-1-4614-1938-9},
  url       = {https://doi.org/10.1007/978-1-4614-1939-6}
}

@article{Moffatt1969,
  author  = {Moffatt, H. K.},
  title   = {The degree of knottedness of tangled vortex lines},
  journal = {Journal of Fluid Mechanics},
  volume  = {35},
  number  = {1},
  pages   = {117--129},
  year    = {1969},
  doi     = {10.1017/S0022112069000991},
  url     = {https://doi.org/10.1017/S0022112069000991}
}

@article{YangEtAl2026,
  author  = {Yang, Qianqian and Gebremichel, Brhane and Ur Rehman, Masood and Yang, Jae Young and Koolen, Jack H.},
  title   = {Sesqui-regular graphs with smallest eigenvalue at least {$-3$}},
  journal = {The Electronic Journal of Combinatorics},
  volume  = {33},
  number  = {3},
  pages   = {P3.20},
  year    = {2026},
  doi     = {10.37236/12174},
  url     = {https://doi.org/10.37236/12174}
}

@article{YangKoolen2025,
  author  = {Yang, Qianqian and Koolen, Jack H.},
  title   = {A structure theory for regular graphs with fixed smallest eigenvalue},
  journal = {Linear Algebra and its Applications},
  volume  = {722},
  pages   = {114--124},
  year    = {2025},
  doi     = {10.1016/j.laa.2025.05.013},
  url     = {https://doi.org/10.1016/j.laa.2025.05.013}
}

@article{Arnold1986,
  author  = {Arnold, Vladimir I.},
  title   = {The asymptotic {Hopf} invariant and its applications},
  journal = {Selecta Mathematica Sovietica},
  volume  = {5},
  number  = {4},
  pages   = {327--345},
  year    = {1986},
  note    = {English translation of the 1974 Russian original}
}

@book{ArnoldKhesin2021,
  author    = {Arnold, Vladimir I. and Khesin, Boris A.},
  title     = {Topological Methods in Hydrodynamics},
  edition   = {2},
  series    = {Applied Mathematical Sciences},
  volume    = {125},
  publisher = {Springer},
  address   = {Cham},
  year      = {2021},
  doi       = {10.1007/978-3-030-74278-2}
}

@article{ChernSimons1974,
  author  = {Chern, Shiing-Shen and Simons, James},
  title   = {Characteristic forms and geometric invariants},
  journal = {Annals of Mathematics},
  volume  = {99},
  number  = {1},
  pages   = {48--69},
  year    = {1974},
  doi     = {10.2307/1971013}
}

@article{Khesin1992,
  author  = {Khesin, Boris A.},
  title   = {Ergodic interpretation of integral hydrodynamic invariants},
  journal = {Journal of Geometry and Physics},
  volume  = {9},
  number  = {1},
  pages   = {101--110},
  year    = {1992},
  doi     = {10.1016/0393-0440(92)90027-X}
}

@article{ArnoldFalkWinther2006,
  author  = {Arnold, Douglas N. and Falk, Richard S. and Winther, Ragnar},
  title   = {Finite element exterior calculus, homological techniques,
             and applications},
  journal = {Acta Numerica},
  volume  = {15},
  pages   = {1--155},
  year    = {2006},
  doi     = {10.1017/S0962492906210018}
}

@article{HuLeeXu2021,
  author  = {Hu, Kaibo and Lee, Young-Ju and Xu, Jinchao},
  title   = {Helicity-conservative finite element discretization
             for incompressible {MHD} systems},
  journal = {Journal of Computational Physics},
  volume  = {436},
  pages   = {110284},
  year    = {2021},
  doi     = {10.1016/j.jcp.2021.110284}
}

@article{HeFarrellHuAndrews2026,
  author  = {He, Mingdong and Farrell, Patrick E. and Hu, Kaibo
             and Andrews, Boris D.},
  title   = {Helicity-preserving finite element discretization
             for magnetic relaxation},
  journal = {SIAM Journal on Scientific Computing},
  volume  = {48},
  number  = {2},
  pages   = {B165--B183},
  year    = {2026},
  doi     = {10.1137/25M1727540}
}

@article{Radhakrishnan1997,
  author  = {Radhakrishnan, Jaikumar},
  title   = {An entropy proof of {Bregman's} theorem},
  journal = {Journal of Combinatorial Theory, Series A},
  volume  = {77},
  number  = {1},
  pages   = {161--164},
  year    = {1997},
  doi     = {10.1006/jcta.1996.2727}
}

@article{SchubelBerwickEvansHirani2024,
  author  = {Schubel, Mark D. and Berwick-Evans, Daniel
             and Hirani, Anil N.},
  title   = {Averaging property of wedge product and naturality
             in discrete exterior calculus},
  journal = {Advances in Computational Mathematics},
  volume  = {50},
  number  = {4},
  pages   = {84},
  year    = {2024},
  doi     = {10.1007/s10444-024-10179-8},
  url     = {https://doi.org/10.1007/s10444-024-10179-8}
}

@book{Greub1978,
  author    = {Greub, Werner},
  title     = {Multilinear Algebra},
  edition   = {2},
  series    = {Universitext},
  publisher = {Springer-Verlag},
  address   = {New York},
  year      = {1978},
  doi       = {10.1007/978-1-4613-9425-9}
}

@book{Weibel1994,
  author    = {Weibel, Charles A.},
  title     = {An Introduction to Homological Algebra},
  series    = {Cambridge Studies in Advanced Mathematics},
  volume    = {38},
  publisher = {Cambridge University Press},
  address   = {Cambridge},
  year      = {1994},
  doi       = {10.1017/CBO9781139644136}
}

@article{Lim2020,
  author  = {Lim, Lek-Heng},
  title   = {{Hodge} Laplacians on Graphs},
  journal = {SIAM Review},
  volume  = {62},
  number  = {3},
  pages   = {685--715},
  year    = {2020},
  doi     = {10.1137/18M1223101}
}

@book{HornJohnson2012,
  author    = {Horn, Roger A. and Johnson, Charles R.},
  title     = {Matrix Analysis},
  edition   = {2},
  publisher = {Cambridge University Press},
  address   = {Cambridge},
  year      = {2012},
  doi       = {10.1017/CBO9781139020411}
}

@book{Diestel2017,
  author    = {Diestel, Reinhard},
  title     = {Graph Theory},
  edition   = {5},
  series    = {Graduate Texts in Mathematics},
  volume    = {173},
  publisher = {Springer},
  address   = {Berlin, Heidelberg},
  year      = {2017},
  doi       = {10.1007/978-3-662-53622-3}
}

\par\bigskip
\begin{samepage}
\footnotesize
\noindent\textsc{Jingqi Li, Fudan University}\par\nopagebreak
\noindent\textit{Email address:}
\href{mailto:23300180152@m.fudan.edu.cn}{\texttt{23300180152@m.fudan.edu.cn}}\par
\end{samepage}
\end{document}